\documentclass [11pt,oneside,a4paper,mathscr]{amsart}

\usepackage{amscd,amsmath,amssymb,euscript}
\usepackage[frame,cmtip,curve,arrow,matrix,line,graph]{xy}

\title{Stability conditions on K3 surfaces}
\author{Tom Bridgeland}

\date{}

\newtheorem{thm}{Theorem}[section]
\newtheorem{conj}[thm]{Conjecture}
\newtheorem{cor}[thm]{Corollary}
\newtheorem{prop}[thm]{Proposition}
\newtheorem{lemma}[thm]{Lemma}
\newenvironment{pf}{\paragraph{Proof}}{\qed\par\medskip}
\theoremstyle{definition}
\newtheorem{defn}[thm]{Definition}

\renewcommand{\leq}{\leqslant}
\renewcommand{\Pi}{\pi}
\renewcommand{\geq}{\geqslant}

\newcommand{\co}{{B}}

\newcommand{\blob}{{\scriptscriptstyle\bullet}}

\newcommand{\Isom}{\operatorname{Aut}}

\newcommand{\Amp}{\operatorname{Amp}}

\newcommand{\End}{\operatorname{End}}

\newcommand{\cl}{\operatorname{c}}
\newcommand{\Pic}{\operatorname{Pic}}

\newcommand{\ha}{\tfrac{1}{2}}
\newcommand{\ei}{\tfrac{1}{8}}
\newcommand{\K}{{{K}}}

\newcommand{\DD}{\operatorname{\mathcal{D}}}

\newcommand{\Ker}{\operatorname{Ker}}
\newcommand{\Coh}{\operatorname{Coh}}

\newcommand{\Aut}{\operatorname{Aut}}
\newcommand{\T}{\operatorname{\mathcal T}}
\newcommand{\isom}{\cong}
\newcommand{\dual}{*}
\newcommand{\td}{\operatorname{td}}
\newcommand{\tensor}{\otimes}
\newcommand{\PP}{\operatorname{\mathbb P}}
\newcommand{\M}{\operatorname{\mathcal M}}
\newcommand{\C}{\mathbb C}

\newcommand{\QQ}{\mathbb Q}

\newcommand{\E}{\mathcal E}
\newcommand{\F}{\mathcal F}

\newcommand{\Z}{\mathbb Z}
\newcommand{\A}{\mathcal A}

\newcommand{\CC}{C}
\newcommand{\OO}{\mathcal O}
\newcommand{\TT}{{T}}

\newcommand{\onto}{\twoheadrightarrow}
\renewcommand{\P}{\mathcal P}

\newcommand{\D}{\operatorname{\DD}}
\newcommand{\Ext}{\operatorname{Ext}}
\newcommand{\Hom}{\operatorname{Hom}}
\newcommand{\eu}{\operatorname{\chi}}
\newcommand{\ch}{\operatorname{ch}}
\renewcommand{\Re}{\operatorname{Re}}
\renewcommand{\Im}{\operatorname{Im}}
\newcommand{\lRa}[1]{\xrightarrow{\ #1\ }}

\newcommand{\lra}{\longrightarrow}

\newcommand{\Stab}{\operatorname{Stab}}
\newcommand{\R}{\mathbb{R}}
\newcommand{\NS}{\operatorname{NS}}
\newcommand{\GL}{\operatorname{GL^+}}
\renewcommand{\S}{S}
\newcommand{\N}{\operatorname{\mathcal{N}}}
\newcommand{\Pol}{\operatorname{\mathcal{P}}}

\newcommand{\grp}{{\tilde{\GL}}(2,\R)}
\newcommand{\m}{m}

\newcommand{\supp}{\operatorname{supp}}
\newcommand{\Q}{\operatorname{\mathcal{Q}}}

\newcommand{\KK}{{\mathcal{K}}}
\newcommand{\U}{{{U}}(X)}
\newcommand{\Ured}{{V}(X)}
\newcommand{\KKo}{{\mathcal{L}}}
\newcommand{\Step}[1]{\smallskip\paragraph{\sc{Step #1}}}

\begin{document}

\begin{abstract}
This paper contains a description of
one connected component of the space of stability
conditions on the bounded derived category of coherent sheaves
on a complex algebraic K3 surface. 
\end{abstract}

\maketitle

\section{Introduction}

The notion of a stability condition on a triangulated category was
introduced in \cite{B} in an effort to understand M. Douglas' work on $\pi$-stability
for D-branes \cite{Do4}.
It was shown in \cite{B} that
to any triangulated category
$\DD$, one can associate a complex manifold
$\Stab(\DD)$ parameterising stability conditions on $\DD$.
Douglas' work suggests that
spaces of stability conditions should be closely related to 
moduli spaces of
superconformal field theories. 
For more on this connection see \cite{B3}.

From a purely
mathematical point of view, spaces of stability conditions are interesting
because they define new invariants of triangulated categories, and
because they
provide a way to introduce geometrical ideas into problems of
homological algebra. Thus for example, the group $\Aut (\DD)$ of
exact autoequivalences of a triangulated category $\DD$ acts on the space
$\Stab(\DD)$ in such a way as to preserve a natural
distance function.

The aim of this paper is to
give a description of one connected component of the space $\Stab(\DD)$ in
the case when $\DD$ is the bounded derived category of coherent sheaves on a complex algebraic K3 surface.
The main result is Theorem \ref{first} below, which
describes a connected component of $\Stab(\DD)$ and shows how its geometry
determines the group $\Aut(\DD)$, which is at present unknown.

\subsection{}
Suppose then that $X$ is an algebraic K3 surface over $\C$.
Following Mukai \cite{Mu2}, one introduces the
extended cohomology lattice of $X$ by using the formula
\[\big((r_1,D_1,s_1),(r_2,D_2,s_2)\big)=D_1\cdot D_2-r_1 s_2-r_2 s_1\]
to define a symmetric bilinear form on the
cohomology ring
\[H^{*}(X,\Z)=H^0(X,\Z)\oplus H^2(X,\Z)\oplus H^4(X,\Z).\]
The resulting lattice $H^*(X,\Z)$ is even and non-degenerate
and has  signature $(4,20)$.
Let $H^{2,0}(X)\subset H^2(X,\C)$
denote the one-dimensional complex subspace spanned by the
class of a nonzero holomorphic two-form $\Omega$ on $X$.
An isometry \[\varphi\colon H^*(X,\Z)\to H^*(X,\Z)\] is called a Hodge
isometry if $\varphi\tensor\C$ preserves this subspace.
The group of Hodge isometries of
$H^*(X,\Z)$ will be denoted $\Aut H^*(X,\Z)$.

Let $\D(X)$ denote the bounded derived category of coherent sheaves on
$X$. The Mukai vector of an object $E\in\D(X)$ is the element of the
sublattice
\[\N(X)=\Z\oplus\NS(X)\oplus\Z=H^*(X,\Z)\cap
\Omega^{\perp}\subset
H^*(X,\C)\]
defined by the formula
\[v(E)=(r(E),c_1(E), s(E))=\ch(E)\sqrt{\td(X)}\in H^*(X,\Z),\]
where $\ch(E)$ is the Chern character of $E$ and $s(E)=\ch_2(E)+r(E)$.

The Mukai bilinear form makes $\N(X)$ into
an even lattice of signature $(2,\rho)$ where
$1\leq\rho\leq 20$ is the Picard number of $X$.
The Riemann-Roch theorem shows that this form is
the negative of the Euler form, that is,
for any pair of objects $E$ and $F$ of $\D(X)$
\[\eu(E,F)=\sum_i (-1)^i \dim_{\C} \Hom_X^i(E,F)=-(v(E),v(F)).\]

A result of Orlov \cite[Proposition 3.5]{Or},
extending work of Mukai \cite[Theorem 4.9]{Mu2}, shows
that every exact autoequivalence of $\D(X)$ induces a Hodge isometry of the
lattice $H^*(X,\Z)$. Thus there is a group homomorphism
\[\varpi\colon \Aut\D(X)\lra\Isom H^*(X,\Z).\]
The kernel of this homomorphism will be denoted $\Aut^0\D(X)$.

Examples of exact
autoequivalences of $\D(X)$ include twists by line bundles and
pullbacks by automorphisms of $X$. A more interesting class of
examples are provided by the twist or reflection functors, first
introduced by Mukai \cite[Proposition 2.25]{Mu2},
and studied in much greater detail
by Seidel and
Thomas \cite{ST}.
Recall that an object
$E\in\D(X)$ is called spherical if
\[\Hom^i_X(E,E)=\begin{cases}
 \,\C &\text{if $i\in\{0,2\}$,} \\
 \, 0 &\text{otherwise.}
\end{cases}\]
Given a spherical object $E\in\D(X),$ the corresponding twist functor
$T_E\in\Aut\D(X)$ is
defined by the triangle
\[ \Hom^{\blob}(E,F)\tensor E\lRa{\alpha} F \lra T_E (F),\]
where the morphism $\alpha:\Hom^{\blob}(E,F)\tensor E\to F$
is the natural evaluation map.

The Riemann-Roch theorem shows that
the Mukai vector of a spherical object lies in the
root system
\[\Delta(X)=\{\delta\in\N(X):(
\delta,\delta)=-2\}.\]
Conversely it is known that for every $\delta\in\Delta(X)$ there exist (infinitely many)
spherical objects with Mukai vector $\delta$. 
Under the homomorphism $\varpi$ the
functor $T_E$ maps to the reflection
\[v\mapsto v+(\delta,v)\delta,\]
where $\delta\in\Delta(X)$ is the Mukai vector of $E$.
Thus the functor
$T_E^{\,2}$ defines an element of $\Aut^0\D(X)$.
Since
$T_E(E)=E[-1]$, this element has infinite order.

\subsection{}
In order to bring geometrical methods to bear on the problem of
computing the group $\Aut \D(X)$ one needs to introduce a space on which it
acts. This is provided by the space of stability conditions on
$\D(X)$.

The precise definition will be recalled in Section \ref{two} below,
but roughly speaking,
a stability condition $\sigma=(Z,\P)$ on a triangulated category
$\DD$ consists of a group
homomorphism
\[Z\colon\K(\DD)\to\C,\]
where $\K(\DD)$ is the Grothendieck
group of $\DD$, and a
collection of
full subcategories $\P(\phi)\subset\DD$, one for each $\phi\in\R$, which
together satisfy a system of axioms.
The map $Z$ is known as the central charge of the stability condition $\sigma$,
and the objects of the subcategory $\P(\phi)$ are said to be semistable of
phase $\phi$ in $\sigma$.

Suppose now that $\DD=\D(X)$ is the derived category of a K3 surface.
A stability condition $\sigma=(Z,\P)$ on $\D(X)$ is said to be
numerical if the central charge $Z$ takes the
form
\[Z(E)=(\Pi(\sigma),v(E))\]
for some vector $\Pi(\sigma)\in\N(X)\tensor\C$.
It was shown in \cite{B} that the set of numerical stability
conditions on $\D(X)$ satisfying a certain technical axiom called
local-finiteness form the points of a complex manifold
$\Stab(X)$. Furthermore the map
\[\Pi\colon\Stab(X)\to\N(X)\tensor\C\]
sending a stability condition to the corresponding vector
$\Pi(\sigma)\in\N(X)\tensor\C$ is continuous.
To understand the image of this map $\Pi$, first define an open subset
\[\Pol(X)\subset\N(X)\tensor\C\]
consisting of those vectors whose real and imaginary parts span positive definite two-planes
in $\N(X)\tensor\R$.
This space has two connected components that are exchanged
by complex conjugation.

Note that $\GL(2,\R)$ acts freely on $\Pol(X)$ by identifying $\N(X)\tensor\C$ with
$\N(X)\tensor\R^2$.
A section of this action is provided by
the submanifold
\[\Q(X)=\{\mho\in \Pol(X):(\mho,\mho)=0,\
(\mho,\bar{\mho})>0,\  r(\mho)=1\}\subset\N(X)\tensor\C,\]
where $r(\mho)$ denotes the projection of $\mho\in\N(X)\tensor\C\subset H^*(X,\C)$ into $H^0(X,\C)$.
The manifold $\Q(X)$ can be identified with the tube domain
\[\big\{\beta+i\omega\in \NS(X)\tensor\C:
\omega^2>0\big\}\]
via the exponential map
\[\mho=\exp(\beta+i\omega)=\big(1,\beta+i\omega,\ha\,(\beta^2-\omega^2) + i(\beta\cdot\omega)\big).\]
Let $\Pol^+(X)\subset\Pol(X)$ denote the connected component containing vectors of the form
$\exp(\beta+i\omega)$ for ample divisor classes
$\omega\in\NS(X)\tensor\R$. For each
$\delta\in\Delta(X)$ let
\[\delta^{\perp}=\{\mho\in\N(X)\tensor\C:(\mho,\delta)=0\}\subset
\N(X)\tensor\C\]
be the
corresponding complex hyperplane. The following is the main result of
this paper.

\begin{thm}
\label{first}
There is a connected component $\Stab^\dagger(X)\subset\Stab(X)$
which is mapped by $\Pi$ onto the open subset
\[\Pol_0^+(X)=\Pol^+(X)\setminus\bigcup_{\delta\in\Delta(X)}\delta^{\perp}
\subset\N(X)\tensor\C.\]
Moreover, the induced map
$\Pi\colon \Stab^\dagger(X)\to\Pol^+_0(X)$
is a covering map, and the subgroup of $\Aut^0 \D(X)$
which preserves the connected component $\Stab^\dagger(X)$ acts freely on
$\Stab^\dagger(X)$ and is the group of deck transformations of $\Pi$.
\end{thm}

The proof of Theorem
\ref{first} depends on a detailed analysis of the open subset
$\U\subset\Stab^\dagger(X)$ consisting of stability conditions in which each
skyscraper sheaf $\OO_x$ is stable of the same phase.
It turns out that these stability conditions can be constructed
explicitly using the method of abstract tilting
first introduced by Happel, Reiten and Smal{\o} \cite{HRS}.

Once one has understood the subset $\U$ it is possible to
show that the components of its boundary
are exchanged by the action of certain autoequivalences of
$\D(X)$. Since these autoequivalences reverse the orientation of the
boundary this quickly leads to a proof of the statement that every
stability condition in $\Stab^\dagger(X)$ is mapped into the closure of
$\U$ by some autoequivalence $\D(X)$.
Theorem \ref{first} then follows relatively easily.

\subsection{}
Unfortunately, Theorem \ref{first} is
not enough to determine the group $\Aut \D(X)$.
The argument of Orlov \cite[Proposition 3.5]{Or}
shows that the image of
the homomorphism
$\varpi$ contains the
index two subgroup \[\Aut^+ H^*(X,\Z)\subset\Aut H^*(X,\Z)\]
consisting of elements which do not exchange the two components of
$\Pol(X)$. But as pointed out by Szendr\H oi \cite{Sz},
it is not known whether the map $\varpi$ is onto.

Note that any element of $\Aut\D(X)$ not mapping into $\Aut^+
H^*(X,\Z)$ could not preserve the component $\Stab^\dagger(X)$.
A related issue is that there may be several components of
$\Stab(X)$ which are permuted by
the action of $\Aut^0\D(X)$. Finally
the space $\Stab^\dagger(X)$ may not be simply-connected.
A proof of the  following geometrical statement
would solve all these problems.

\begin{conj}
\label{conj}
The action of $\Aut \D(X)$ on $\Stab(X)$ preserves the connected
component $\Stab^\dagger(X)$. Moreover $\Stab^\dagger(X)$ is simply-connected.
Thus there is a short exact
sequence of groups
\[1\lra \pi_1 \Pol_0^+(X)\lra
\Aut \D(X)\lRa{\varpi} \Aut^+ H^*(X,\Z)\lra
1.\]
\end{conj}

It seems reasonable to hope that a more detailed analysis of
$\Stab(X)$, and in particular
its natural distance function, will lead to a proof of this conjecture.

\subsection*{Acknowledgements} 
Thanks first of all to Mike Douglas whose work on $\Pi$-stability
provided the basic idea for this paper. I'm also very grateful to So Okada, Bal\'azs
Szendr\H oi, Hokotu Uehara and particularly Daniel Huybrechts
for pointing out various gaps and errors in  earlier versions.

% ***************************************************************************
% ***************************************************************************
% ***************************************************************************
% ***************************************************************************

\section{Stability conditions}
\label{two}

The notion of a stability condition on a triangulated category was
introduced in \cite{B}. The next three sections contain a summary of
the contents of that paper together with
a few additional statements that will be needed later on.

\begin{defn}
\label{pemb}
A stability condition $\sigma=(Z,\P)$ on a triangulated category $\DD$
consists of
a linear map
$Z\colon\K(\DD)\to\C$ called the {central charge},
and full additive
subcategories $\P(\phi)\subset\DD$ for each $\phi\in\R$,
satisfying the following axioms:
\begin{itemize}
\item[(a)] if $0\neq E\in \P(\phi)$ then $Z(E)=\m(E)\exp(i\pi\phi)$ for some
 $\m(E)\in\R_{>0}$,
\item[(b)] for all $\phi\in\R$, $\P(\phi+1)=\P(\phi)[1]$,
\item[(c)] if $\phi_1>\phi_2$ and $A_j\in\P(\phi_j)$ then $\Hom_{\DD}(A_1,A_2)=0$,
\item[(d)] for $0\neq E\in\DD$ there is a finite sequence of real
numbers
\[\phi_1>\phi_2> \cdots >\phi_n\]
and a collection of triangles
\[
\xymatrix@C=.5em{
0_{\ } \ar@{=}[r] & E_0 \ar[rrrr] &&&& E_1 \ar[rrrr] \ar[dll] &&&& E_2
\ar[rr] \ar[dll] && \ldots \ar[rr] && E_{n-1}
\ar[rrrr] &&&& E_n \ar[dll] \ar@{=}[r] &  E_{\ } \\
&&& A_1 \ar@{-->}[ull] &&&& A_2 \ar@{-->}[ull] &&&&&&&& A_n \ar@{-->}[ull] 
}
\]
with $A_j\in\P(\phi_j)$ for all $j$.

\end{itemize}
\end{defn}

Given a stability condition as in the definition, each subcategory
$\P(\phi)$ is abelian. The nonzero objects of $\P(\phi)$ are said to be
{semistable of phase $\phi$ in $\sigma$}, and the simple objects
of $\P(\phi)$ are said to be stable.
It is an easy exercise to check that the decompositions of
a nonzero object $0\neq E\in\DD$ given by
axiom (d) are uniquely defined up to isomorphism.
The objects $A_j$ will be called the {semistable factors} of $E$
with respect to $\sigma$.
Write
$\phi^+_{\sigma}(E)=\phi_1$ and
$\phi^-_{\sigma}(E)=\phi_n$; clearly
$\phi^-_{\sigma}(E)\leq\phi^+_{\sigma}(E)$ with equality
holding precisely when $E$ is semistable in $\sigma$.
The {mass} of $E$ is defined to be
the positive real number
$m_{\sigma}(E)=\sum_i |Z(A_i)|$.
By the triangle inequality
one has $m_{\sigma}(E)\geq |Z(E)|$.
When the stability condition $\sigma$ is clear from the context
I often
drop it from the
notation and write $\phi^{\pm}(E)$ and $m(E)$.

For any interval $I\subset\R$, define $\P(I)$ to be the extension-closed
subcategory of $\DD$ generated by the subcategories $\P(\phi)$ for
$\phi\in I$. Thus, for example, the full subcategory $\P((a,b))$
consists of the zero objects of $\DD$ together with those
objects $0\neq E\in\DD$ which satisfy
$a<\phi^-(E)\leq\phi^+(E)<b$.

To prove nice results it is necessary to put one extra condition on
stability conditions. A stability condition is called {locally finite} if
there is some $\epsilon>0$ such that each quasi-abelian category
$\P((\phi-\epsilon,\phi+\epsilon))$ is of finite length. For more
details on this see
\cite{B}. 
If $\sigma=(Z,\P)$ is locally finite then each subcategory
$\P(\phi)$ has finite length, so that every semistable object has a
finite Jordan-H{\"o}lder filtration into stable factors of the same
phase.

The set $\Stab(\DD)$ of locally-finite stability conditions on a fixed
triangulated category $\DD$ has a natural topology induced by the
generalised metric
\[d(\sigma_1,\sigma_2)=\sup_{0\neq E\in\DD}
\bigg\{|\phi^-_{\sigma_2}(E)-\phi^-_{\sigma_1}(E)|,|\phi^+_{\sigma_2}(E)
-\phi^+_{\sigma_1}(E)|
,|\log \frac{m_{\sigma_2}(E)}{m_{\sigma_1}(E)}|\bigg\}
\in[0,\infty].\]
With this topology the functions $m_{\sigma}(E)$ and
$\phi^{\pm}_{\sigma}(E)$ are continuous for every nonzero object
$E\in\DD$.
It follows that the subset of $\Stab(\DD)$ where a given
object $E\in\DD$ is semistable is a closed subset.

\begin{lemma}\cite[Lemma 8.2]{B}
\label{groupactions}
The generalised metric space $\Stab(\DD)$ carries a right action of the
group
$\grp$, the universal
cover of $\GL(2,\R)$, and a left action by isometries of the 
group
$\Aut(\DD)$ of exact autoequivalences of $\DD$. These two actions
commute.
\end{lemma}

\begin{pf}
First note that the group $\grp$ can be thought of as the
set of pairs $(T,f)$ where $f\colon \R\to\R$ is an
increasing map with $f(\phi+1)=f(\phi)+1$, and  $T\colon \R^2\to\R^2$
is an orientation-preserving linear isomorphism, such that the induced maps on
$S^1=\R/2\Z=(\R^2\setminus\{0\})/\R_{>0}$ are the same.

Given a stability condition $\sigma=(Z,\P)\in\Stab(\DD)$, and a
pair $(T,f)\in\grp$, define a new stability
condition $\sigma'=(Z',\P')$ by setting $Z'=T^{-1}\circ Z$ and
$\P'(\phi)=\P(f(\phi))$. Note that the semistable objects of the
stability conditions $\sigma$
and $\sigma'$ are the same, but the phases have been
relabelled.

For the second action, note that an element
$\Phi\in\Aut(\DD)$ induces an automorphism $\phi$ of $\K(\DD)$. If
$\sigma=(Z,\P)$ is a stability condition on $\DD$, define $\Phi(\sigma)$
to be the stability condition $(Z\circ\phi^{-1},\P')$, where
$\P'(t)$=$\Phi(\P(t))$. The reader can easily check that this action
is by isometries and commutes with the first.
\end{pf}

If $\sigma$ and $\tau$ are stability conditions
on a triangulated category $\DD$ define
\[f(\sigma,\tau)=\sup_{0\neq E\in\DD}
\bigg\{|\phi^-_{\tau}(E)-\phi^-_{\sigma}(E)|,|\phi^+_{\tau}(E)
-\phi^+_{\sigma}(E)|\bigg\}\in[0,\infty].\]
For any connected component $\Stab^*(X)\subset\Stab(\DD)$ the function
$f\colon\Stab^*(X)\times\Stab^*(X)\to\R$ is continuous and finite.
It appeared in \cite{B} as a generalised metric on the space of
slicings.

\begin{lemma}\cite[Lemma 6.4]{B}
\label{rachael}
If $\sigma$ and $\tau$ are stability conditions
on a triangulated category
$\DD$ with the same central charge, and $f(\sigma,\tau)< 1$,
then $\sigma=\tau$.
\qed
\end{lemma}

The main result of \cite{B} is the following deformation result. The idea
is that if $\sigma=(Z,\P)$ is a stability condition on a triangulated
category $\DD$ and one
deforms $Z$ to a new group homomorphism $W\colon\K(\DD)\to\C$ in such a way that
the phase of each stable object in $\sigma$ changes in a uniformly
bounded way, then it is possible to  define new classes of semistable objects
$\Q(\psi)\subset\DD$ so that $(W,\Q)$ is a stability condition on $\DD$.

\begin{thm}\cite[Theorem 7.1]{B}
\label{biggy}
Let $\sigma=(Z,\P)$ be a locally-finite
stability condition on a triangulated category $\DD$. Take $0<\epsilon<\ei$ such
that each of the quasi-abelian categories
$\P((t-4\epsilon,t+4\epsilon))\subset\D$ is of finite length.
If 
$W\colon\K(\DD)\to\C$ is a group homomorphism satisfying
\[|W(E)-Z(E)|<\sin(\pi\epsilon)|Z(E)|\]
for all objects $E\in\DD$ which are stable in $\sigma$,
then there is a locally-finite stability condition
$\tau=(W,\Q)$ on $\DD$ with
$f(\sigma,\tau)<\epsilon$.
\qed
\end{thm}

In fact in \cite[Theorem 7.1]{B} the inequality was assumed for all objects $E\in\DD$ which are semistable in $\sigma$, but
it is clear that it is enough to check it for stable objects, since any semistable object has a finite filtration by stable objects of the same phase.

% ***************************************************************************
% ***************************************************************************
% ***************************************************************************
% ***************************************************************************

\section{T-structures and stability functions}
A stability condition on a triangulated category $\DD$ consists of a
t-structure on $\DD$ together with a stability function on its heart.
This statement is made precise in
Proposition \ref{pg} below.
Readers unfamiliar with the concept of a
t-structure should consult \cite{BBD,GM}.
The following easy result is a good exercise.

\begin{lemma}
\label{ll}
A bounded t-structure is determined by its heart.
Moreover, if  $\A\subset\DD$ is a full additive
subcategory of a triangulated category
$\DD$,
then $\A$ is the heart of a bounded t-structure on $\DD$ if and
only if the
following two conditions hold:
\begin{itemize}
\item[(a)]if $A$ and $B$ are objects of
$\A$ then $\Hom_{\DD}(A,B[k])=0$ for $k<0$,

\item[(b)]
for every nonzero object $E\in\DD$ there are integers $m<n$
and a collection of triangles
\[
\xymatrix@C=.2em{
0_{\ } \ar@{=}[r] & E_{m} \ar[rrrr] &&&& E_{m+1} \ar[rrrr] \ar[dll] &&&& E_{m+2}
\ar[rr] \ar[dll] && \ldots \ar[rr] && E_{n-1}
\ar[rrrr] &&&& E_n \ar[dll] \ar@{=}[r] &  E_{\ } \\
&&& A_{m+1} \ar@{-->}[ull] &&&& A_{m+2} \ar@{-->}[ull] &&&&&&&& A_n \ar@{-->}[ull] 
}
\]
with $A_i[i]\in\A$ for all $i$.\qed
\end{itemize}
\end{lemma}

In analogy with the standard t-structure on the derived category of an
abelian category, the objects $A_i[i]$ of $\A$
are called the cohomology objects
of $A$ in the given t-structure, and one often writes $H^i(E)=A_i[i]$.

A very useful method for constructing t-structures is provided
by the idea of
tilting with respect to a torsion pair, first introduced by
D. Happel, I. Reiten and
S. Smal{\o} \cite{HRS}.

\begin{defn}
\label{tors}
A torsion pair in an abelian category
$\A$ is a pair of full subcategories $(\T,\F)$ of $\A$
which satisfy
$\Hom_{\A}(T,F)=0$ for $T\in \T$ and $F\in\F$,
and such that every object $E\in\A$ fits into a  short exact sequence
\[0\lra T\lra E\lra F\lra 0\] for some pair of objects $T\in\T$ and
$F\in  \F$.
\end{defn}

The objects of $\T$ and $\F$ are called torsion and torsion-free
respectively. The following Lemma is pretty-much immediate from
Lemma \ref{ll}.

\begin{lemma}\cite[Proposition 2.1]{HRS}
Suppose $\A$ is the heart of a bounded
t-structure on a
triangulated category $\DD$. Given an object $E\in\DD$ let $H^i(E)\in\A$
denote the
$i$th cohomology object of $E$ with respect to this t-structure.
Suppose $(\T,\F)$ is a torsion pair in $\A$. Then 
the full subcategory
\[\A^{\sharp}=\big\{E\in \DD:H^i(E)=0\text{ for }i\notin\{-1,0\},
H^{-1}(E)\in\F\text{ and } H^{0}(E)\in\T\big\}\]
is the heart of a bounded t-structure on $\DD$.
\qed
\end{lemma}

One says that $\A^{\sharp}$ is obtained from the category $\A$ by tilting with
respect to the torsion pair $(\T,\F)$. Note that the pair
$(\F[1],\T)$
is a torsion pair in $\A^{\sharp}$
and that tilting with respect to this pair
gives back the original category $\A$ with a shift.

\begin{defn}
A stability function on
an abelian  category $\A$ is a group homomorphism $Z\colon \K(\A)\to\C$
such that
\[0\neq E\in\A \implies Z(E)\in\R_{>0}\,\exp({i\pi\phi(E)})\text{ with
}0<\phi(E)\leq 1.\]
The real number $\phi(E)\in(0,1]$ is called the phase of the object $E$.
\end{defn}

Such stability functions were called centered in \cite{B} but in this paper we
shall make no use of non-centered stability functions.
A nonzero object $E\in\A$ is said to be
semistable
with respect to a stability function $Z$
if
\[0\neq A\subset E\implies  \phi(A)\leq\phi(E).\]
The stability function $Z$ is said to have the Harder-Narasimhan property if
every nonzero object $E\in\A$ has
a finite filtration
\[0=E_0\subset E_1\subset \cdots\subset E_{n-1}\subset E_n=E\]
whose factors $F_j=E_j/E_{j-1}$ are semistable objects of $\A$ with
\[\phi(F_1)>\phi(F_2)>\cdots>\phi(F_n).\]
A simple sufficient condition for the existence of Harder-Narasimhan
filtrations was given in \cite[Proposition 2.4]{B}.

\begin{prop}\cite[Proposition 5.3]{B}
\label{pg}
To give a stability condition on a triangulated category $\DD$ is
equivalent to giving a bounded t-structure on $\DD$ and a 
stability function on its heart which
has the Harder-Narasimhan property.
\end{prop}

\begin{pf}
If $\sigma=(Z,\P)$ is a stability condition on
$\DD$, the abelian subcategory  $\A=\P((0,1])$
is the heart of a bounded t-structure on $\DD$.
The central charge $Z$ defines
a  stability function on $\A$ and it is easy to check
that the corresponding semistable objects are
precisely the nonzero objects of the
subcategories $\P(\phi)$ for $0<\phi\leq 1$.
The decompositions of objects of $\A$ given by Definition
\ref{pemb}(d) are Harder-Narasimhan filtrations.

For the converse, suppose $\A$ is the heart of a bounded t-structure
on $\DD$, and $Z\colon\K(\A)\to\C$ is a  stability function on $\A$
with the Harder-Narasimhan property. Define a stability condition
$\sigma=(Z,\P)$ on $\DD$ as follows.
For each $\phi\in(0,1]$ let $\P(\phi)$ be the
full additive subcategory of $\DD$
consisting of
semistable objects
of $\A$ with phase $\phi$, together with the zero objects of $\DD$. 
Condition (b) of Definition \ref{pemb} then
determines $\P(\phi)$ for all $\phi\in\R$ and conditions (a) and (c)
are easily
verified. Given a nonzero object $E\in\DD$, a filtration as in
Definition \ref{pemb}(d)
can be obtained by combining the decompositions of Lemma \ref{ll} with the
Harder-Narasimhan filtrations of the cohomology objects of $E$.
\end{pf}

In Section \ref{con} this result will be combined with the
method of tilting and used to construct examples of stability conditions on
K3 surfaces.

% ***************************************************************************
% ***************************************************************************
% ***************************************************************************
% ***************************************************************************

\section{Stability conditions on varieties}

In this paper the triangulated category $\DD$ will always be the
bounded derived category of coherent sheaves $\D(X)$ on a smooth projective
variety $X$ over the complex numbers. Any such
category is of finite type, that is, for any pair of
objects $E$ and $F$  in $\DD$ the morphism space
$\bigoplus_i \Hom_{X}(E,F[i])$
is a finite-dimensional vector space over $\C$. The Euler characteristic
\[\eu(E,F)=\sum_i (-1)^i \dim_{\C} \Hom_{X}(E,F[i]),\]
then defines a bilinear form on the Grothendieck group $\K(X)$.
Serre duality shows that the left- and right-radicals
$\K(X)^{\perp}$ and $\smash{^\perp}{\K(X)}$ are the same, so that the
Euler form descends to a nondegenerate form $\eu(-,-)$ on the
numerical Grothendieck group 
\[\N(X)=\K(X)/\K(X)^{\perp}.\] The Riemann-Roch theorem shows that
this free abelian group has finite rank. A stability condition
$\sigma=(Z,\P)$ is said to be
{numerical} if the central
charge $Z\colon\K(X)\to\C$ factors through the quotient group $\N(X)$.
An equivalent condition is that the central charge $Z$ takes
the form
\[Z(E)=-\eu(\Pi(\sigma),E)\]
for some vector
$\Pi(\sigma)\in\N(X)\tensor\C$.
Here $E$ is simultaneously representing an element of
$\K(X)$ and the quotient $\N(X)$. The sign is just to
ensure agreement with later conventions.

The set of all
locally finite numerical stability conditions on $\D(X)$
with its natural topology will be denoted $\Stab(X)$. The group $\Aut \D(X)$ of exact autoequivalences
of $\D(X)$ acts on $\Stab(X)$. 
Theorem \ref{biggy} leads to the following result.

\begin{thm} \cite[Cor. 1.3]{B}
\label{cath}
For each connected component $\Stab^*(X)\subset \Stab(X)$ there is a
linear subspace $V\subset \N(X)\tensor\C$ such that
\[\pi\colon\Stab^*(X)\lra\N(X)\tensor\C\]
is a local
homeomorphism onto an open subset of the subspace $V$.
In particular $\Stab^*(X)$ is a finite-dimensional
complex manifold.
\end{thm}

In all known examples, the subspace $V$ of Theorem
\ref{cath} is actually equal to $\N(X)\tensor\C$. Since it is not known whether this is always the case
the following definition will be useful.

\begin{defn}
\label{full}
A connected component $\Stab^*(X)\subset\Stab(X)$ will be called full if the subspace $V$ of Theorem \ref{cath}
is equal to $\N(X)\tensor\C$. A stability condition $\sigma\in\Stab(X)$ is full if it lies in a full component.
\end{defn}

Note that if a stability condition $\sigma\in\Stab(X)$ is full and $\Phi\in\Aut\D(X)$ is an autoequivalence then the stability condition $\Phi(\sigma)$ is also
full.

Later on it will be important to be able to choose the constant $\epsilon$ appearing in Theorem \ref{biggy} uniformly.
This can be done for full stability conditions.
To do so one first considers discrete stability conditions.

\begin{defn}
A stability condition $\sigma=(Z,\P)$ on $\D$ is called discrete if the image of $Z\colon\K(\D)\to\C$ is a discrete subgroup.
\end{defn}

\begin{lemma}
\label{disc}
Suppose $\sigma=(Z,\P)$ is a  discrete stability condition and fix $0<\epsilon<\ha$. Then for each $\phi\in\R$ the quasi-abelian 
category $\P((\phi-\epsilon,\phi+\epsilon))$ is of finite length. In particular $\sigma$ is locally finite.
\end{lemma}

\begin{pf}
Fix $\phi\in\R$ and set $\A=\P((\phi-\epsilon,\phi+\epsilon))$. The central charge of any nonzero object $A\in\A$ lies in the sector
\[S=\{z=r\exp(i\pi\psi):r>0 \text{ and }\phi-\epsilon<\psi<\phi+\epsilon\}\]
which is strictly smaller than a half-plane in $\C$. Set $f(A)=\Re \big(\exp(-i\pi\phi)Z(A)\big)$. Then $f(A)>0$ for all nonzero objects
$A\in\A$. If
\[0\lra A\lra B\lra C\lra 0\]
is a strict short exact sequence in $\A$ then $f(B)=f(A)+f(C)$. Thus for a given object $E\in\A$ the central charges of
all sub and quotient objects lie in the bounded region
\[\{z\in S:\Re \big(\exp(-i\pi\phi)z\big)< f(E)\}\]
Since $\sigma$ is discrete there are only finitely many possibilities, and so any chains of sub or quotient objects must terminate.
\end{pf}

\begin{lemma}
\label{fully}
Suppose $\sigma\in\Stab(X)$ is a full stability condition and fix $0<\epsilon<\ha$. Then for each $\phi\in\R$ the quasi-abelian 
category $\P((\phi-\epsilon,\phi+\epsilon))$ is of finite length. 
\end{lemma}

\begin{pf}
The basic point is that there exist discrete stability conditions
in $\Stab(X)$ arbitrarily close to $\sigma=(Z,\P)$. Indeed, one can approximate $\pi(\sigma)\in\N(X)\tensor\C$ arbitrarily closely by
points in $\N(X)\tensor\QQ[i]$. By the fullness assumption these points can be lifted to stability conditions $\tau=(W,\Q)$ lying arbitrarily close to $\sigma$. The resulting stability conditions will clearly be discrete.
But by \cite[Lemma 6.1]{B}, if $f(\sigma,\tau)<\eta$ then $\P((\phi-\epsilon,\phi+\epsilon))\subset \Q((\phi-\epsilon-\eta,\phi+\epsilon+\eta))$ so the result follows from Lemma \ref{disc}.
\end{pf}

% ***************************************************************************
% ***************************************************************************
% ***************************************************************************
% ***************************************************************************

\section{Sheaves on K3 surfaces}

This section contains some basic results about coherent sheaves on K3
surfaces. These are mostly taken from Mukai's excellent paper
\cite{Mu2}.
Throughout notation will be as in the introduction. In
particular $X$ is a fixed algebraic K3 surface over $\C$ and $\D(X)$
is the corresponding bounded derived category of coherent sheaves. 

There are various notions of stability for sheaves on a K3 surface,
the most basic of which is slope-stability.
Recall that if $\omega\in\NS(X)$ is an ample divisor class, one defines the
slope $\mu_\omega(E)$ of a torsion-free
sheaf $E$ on $X$ to be the quotient
\[\mu_\omega(E)=\frac{\cl_1(E)\cdot \omega}{r(E)}.\]
A torsion-free sheaf $E$ is said to be $\mu_\omega$-semistable if
$\mu_\omega(A)\leq\mu_\omega(E)$ for every subsheaf
$0\neq A\subset E$. If
the inequality is always strict when $A$ is of strictly smaller rank then $E$ is said to be $\mu_\omega$-stable.
This definition
can be extended to include any class $\omega$ in the ample
cone of $X$,
\[\Amp(X)=\{\omega\in \NS(X)\tensor \R : \omega^2>0\text{ and }\omega\cdot
C>0\text{ for any curve }C\subset X\}.\]

One very important point to realise is
that slope-stability does not arise from a stability
condition on $\D(X)$. The function $Z(E)=-\cl_1(E)\cdot \omega +
ir(E)$ is not a stability function on the category of coherent sheaves on
$X$, because it is zero on any sheaf supported in dimension zero. Thus
constructing examples of stability conditions on $\D(X)$ is a
non-trivial problem, which will be tackled in the next section.

Sending an object $E\in\D(X)$ to its Mukai vector
\[v(E)=(r(E),c_1(E), s(E))=\ch(E)\sqrt{\td(X)}\in H^*(X,\Z)\]
identifies the
numerical Grothendieck group of $X$ with the sublattice
\[\N(X)=\Z\oplus\NS(X)\oplus\Z=H^*(X,\Z)\cap
\Omega^{\perp}\subset
H^*(X,\C)\]
defined in the introduction.
The Riemann-Roch theorem then shows that the Mukai bilinear form
induces the negative of the Euler form on $\N(X)$
\cite[Proposition
2.2]{Mu2}, so that
for any pair of objects $E$ and $F$ of $\D(X)$
\[\eu(E,F)=\sum_i (-1)^i \dim_{\C} \Hom_X^i(E,F)=-(v(E),v(F)).\]
Recall also \cite[Proposition
2.3]{Mu2} that Serre duality gives isomorphisms
\[ \Hom^i_X(E,F)\isom \Hom^{2-i}_X(F,E)^{\dual}. \]

If the objects $E$ and $F$ lie in the heart of some t-structure on
$\D(X)$ these spaces vanish for
$i<0$ and hence also for $i>2$. This is the case for example
when $E$ and $F$ are both
sheaves, or if $E$ and $F$ are semistable of the same phase in some
stability condition. 
In this situation, combining Riemann-Roch and Serre duality
often allows one to determine the spaces
$\Hom^i_X(E,F)$. The next Lemma is a good example.

\begin{lemma}\cite[Corollary 2.5]{Mu2}
\label{cleave}
If $E\in\D(X)$ is stable in some
stability condition on $X$, or is a $\mu_{\omega}$-stable sheaf for some
$\omega\in\Amp(X)$, then
\[\dim_{\C} \Hom^1_X(E,E) =2+v(E)^2\geq 0,\]
with equality precisely when $E$ is spherical.
\end{lemma}

\begin{pf}
Any stable object satisfies $\Hom_X(E,E)=\C$, so that by Serre
duality, one has $\Hom^2_X(E,E)=\C$ also, and the given inequality
follows from Riemann-Roch.
\end{pf}

The following Lemma is essentially
due to Mukai.

\begin{lemma}\cite[Corollary 2.8]{Mu2}
\label{mukai}
Suppose $\A\subset\D(X)$ is the heart of a bounded t-structure and
\[0\lra A\lra B\lra C\lra 0\]
is a short exact sequence in $\A$ with $\Hom_X(A,C)=0$. Then
\[\dim_\C \Hom^1_X(A,A)+\dim_\C \Hom^1_X(C,C)\leq\dim_\C
\Hom^1_X(B,B).\]
\end{lemma}

\begin{pf}
Given objects $E$ and $F$ of $\A$ write
\[(E,F)^i=\dim_\C \Hom^i_X(E,F).\]
Thus $(E,F)^i=0$ unless $0\leq i\leq
2$, and Serre duality gives $(E,F)^i=(F,E)^{2-i}$.
The stated inequality is equivalent to
\[(A,A)^0+(C,C)^0+(C,A)^0\leq (C,A)^1+(B,B)^0,\]
which follows from the existence of an exact sequence
\[0\to \Hom_X(C,A)\to \End_X(B)\to \End_X(A)\oplus\End_X(C)\to
\Hom^1_X(C,A).\]
To see where this sequence comes from note that because $(A,C)^0=0$,
any endomorphism of $B$ induces an endomorphism of the whole
triangle $A\to B\to C$, which is just a pair of endomorphisms of
$A$ and $C$ preserving the class of the connecting morphism $C\to A[1]$.
\end{pf}

The following important result of Yoshioka
gives existence of semistable sheaves.

\begin{thm}{\rm (Yoshioka)}
\label{yosh}
Suppose $v\in\N(X)$ satisfies $v^2\geq -2$ and $r(v)>0$. Then for any
ample divisor class $\omega\in\Amp(X)$, there
are torsion-free $\mu_{\omega}$-semistable sheaves on $X$ with
Mukai vector $v$.
\end{thm}

\begin{pf}
Since one is looking only for semistable sheaves it is enough to check
the case when $v\in\N(X)$ is primitive. Using the wall and chamber
structure of $\Amp(X)$ (see for example \cite{Yo2})
one can also deform $\omega$ a little
so that it is general for $v$, in the sense that all $\mu_{\omega}$-semistable
sheaves with Mukai vector $v$ are actually $\mu_{\omega}$-stable. The
result then follows from \cite[Theorem 8.1]{Yo1}.
\end{pf}

% ***************************************************************************
% ***************************************************************************
% ***************************************************************************
% ***************************************************************************

\section{Constructing stability conditions}
\label{con}

Having dealt with the preliminaries it is now possible to make a start
on the proof of Theorem \ref{first}. For the next nine sections $X$
will be a fixed algebraic K3 surface over $\C$.
In the next two sections we use the method of tilting to
construct examples of stability conditions on the bounded derived category
of coherent sheaves $\D(X)$.

Take a pair of $\R$-divisors $\beta,\omega\in \NS(X)\tensor\R$
such that $\omega\in\Amp(X)$ lies in the ample cone.
As in the introduction, define a group homomorphism $Z\colon\N(X)\to\C$ by
the formula
\[Z(E)=(\exp(\beta+i\omega),v(E)).\]
A little rewriting shows that if $E\in\D(X)$ has
Mukai vector $(r,\Delta,s)$ with $r\neq 0$ then
\begin{equation*}
\tag{$\star$}Z(E)=\frac{1}{2r}\bigg((\Delta^2-2rs)+r^2\omega^2-(\Delta-r\beta)^2\bigg)
+i(\Delta-r\beta)\cdot\omega
\end{equation*}
which reduces to
$Z(E)=(\Delta\cdot \beta-s) +i(\Delta\cdot \omega)$ when $r=0$.

Every torsion-free sheaf $E$ on $X$
has a unique Harder-Narasimhan filtration
\[0=E_0\subset E_1\subset \cdots \subset E_{n-1}\subset E_n=E\]
whose factors $F_i=E_i/E_{i-1}$ are $\mu_\omega$-semistable torsion-free
sheaves with descending slope $\mu_\omega$.
Truncating Harder-Narasimhan filtrations at the point
$\mu_{\omega}=\beta\cdot\omega$ leads to the following statement.

\begin{lemma}
\label{mckay}
For any pair $\beta,\omega\in\NS(X)\tensor\R$ with $\omega\in\Amp(X)$
there is a unique torsion
pair $(\T,\F)$ on the category $\Coh(X)$ such that
$\T$ consists of sheaves whose torsion-free parts have
 $\mu_{\omega}$-semistable
Harder-Narasimhan factors of slope $\mu_{\omega}>\beta\cdot
\omega$ and $\F$ consists of
torsion-free
sheaves on $X$ all of
whose $\mu_{\omega}$-semistable
Harder-Narasimhan factors have slope $\mu_{\omega}\leq\beta\cdot
\omega$.
\qed
\end{lemma}

Tilting with respect to the torsion pair of Lemma \ref{mckay}
gives a bounded t-structure on $\D(X)$ with
heart
\[\A(\beta,\omega)=\big\{E\in \D(X):H^i(E)=0\text{ for }i\notin\{-1,0\},
H^{-1}(E)\in\F\text{ and } H^{0}(E)\in\T\big\}.\]
Note that $\A(\beta,\omega)$ does not really depend on $\beta$, only on $\beta\cdot\omega$.
Note also that all torsion sheaves on $X$ are objects of $\T$ and hence also
of the abelian category $\A(\beta,\omega)$.

\begin{lemma}
\label{slopefun}
Take a pair $\beta,\omega\in \NS(X)\tensor\R$ with
$\omega\in\Amp(X)$. Then the group homomorphism $Z$
defined above is a
stability function on the abelian category $\A(\beta,\omega)$
providing $\beta$ and $\omega$ are chosen so that
for all spherical sheaves $E$ on $X$ one has $Z(E)\notin \R_{\leq
0}$. This holds in particular whenever $\omega^2>2$.
\end{lemma}

\begin{pf}
It is clear that $Z(E)$ lies in the upper half-plane
for every sheaf supported on a curve, and every torsion-free
semistable sheaf $E$ with $\mu_\omega(E)>\beta\cdot
\omega$. Moreover,
if $E$ is supported in dimension zero then
$Z(E)\in\R_{<0}$. Similarly, if $E$ is torsion-free
with $\mu_\omega(E)<\beta\cdot
\omega$ then $Z(E[1])$ lies in the upper half-plane.

The only
non-trivial part is to check that if a torsion-free
$\mu_{\omega}$-semistable sheaf $E$ satisfies
$(\Delta-r\beta)\cdot \omega=0$ then $Z(E)\in\R_{>0}$.
It is enough to check
this when $E$ is $\mu_{\omega}$-stable.
By Lemma \ref{cleave}, the Mukai vector $v(E)=(r,\Delta,s)$ satisfies
\[\Delta^2-2rs=v(E)^2\geq -2\]
with equality precisely when $E$ is spherical.
Since $(\Delta-r\beta)\cdot \omega=0$
the Hodge index theorem gives $(\Delta-r\beta)^2\leq 0$,
so if $E$ is not spherical then $Z(E)$ lies on the
positive real axis. If  $\omega^2>2$
then $Z(E)$ lies on the positive real axis even if $E$ is spherical.
\end{pf}

Suppose the pair $\beta,\omega\in \NS(X)\tensor\R$
satisfy the condition of Lemma \ref{slopefun}.
If the corresponding stability function $Z$
has the Harder-Narasimhan property,
Proposition  \ref{pg} shows that there is a uniquely-defined
stability condition $\sigma=(Z,\P)$ on $\D(X)$ with heart
$\P((0,1])=\A(\beta,\omega)$ and central charge $Z$.
In Section \ref{U} it will be shown that in fact this Harder-Narasimhan property always holds.
In the next section a weaker result will be proved, namely that the Harder-Narasimhan property 
holds when $\beta$ and
$\omega$ are rational. 

The following observation will be used later
to characterise the potential stability conditions constructed above.

\begin{lemma}
\label{balloon}
For any pair $\beta,\omega\in\NS(X)\tensor\R$ with $\omega\in\Amp(X)$,
and any point $x\in X$,
the skyscraper sheaf $\OO_x$ is a simple object of the abelian
category $\A(\beta,\omega)$.
\end{lemma}

\begin{pf}
As remarked above, any torsion sheaf on $X$ lies in the torsion
subcategory $\T$ and hence in
the abelian category $\A(\beta,\omega)$. Suppose
\[0\lra A\lra \OO_x\lra B\lra 0\]
is a short exact sequence in $\A(\beta,\omega)$.
Taking cohomology gives an exact
sequence of sheaves on $X$
\[0\lra H^{-1}(B)\lra H^0(A)\lra \OO_x\lra H^0(B)\lra 0.\]
Note that $H^{-1}(B)$ is torsion-free.
It follows that the $\mu_{\omega}$-semistable
factors of $H^{-1}(B)$ and
$H^0(A)$ have the same slope. This contradicts the definition of the
category $\A(\beta,\omega)$
unless $H^{-1}(B)=0$, in which case either $A$ or $B$
must be zero.
\end{pf}

% ***************************************************************************
% ***************************************************************************
% ***************************************************************************
% ***************************************************************************

\section{The Harder-Narasimhan property}

This section fills a gap in the original version of this paper. I am very grateful to D. Huybrechts for
pointing out the error, and to E. Macri and P. Stellari who independently obtained Proposition \ref{stella} below.

Take $\beta,\omega\in\NS(X)\tensor\R$ with $\omega\in\Amp(X)$ and let $Z\colon \N(X)\to\C$ be the group homomorphism
\[Z(E)=(\exp(\beta+i\omega),v(E))\]
defined in the previous section. Suppose that for all spherical sheaves $E$ on $X$ one has
$Z(E)\notin \R_{\leq 0}$. According to Lemma \ref{slopefun} the map $Z$ then defines a stability function
on the abelian category $\A(\beta,\omega)\subset\D(X)$.

\begin{prop}
\label{stella}
If $\beta,\omega\in\NS(X)\tensor\QQ$ are rational then the stability function $Z$ on $\A(\beta,\omega)$ has the Harder-Narasimhan property.
The resulting stability conditions on $\D(X)$ are locally-finite.
\end{prop}

\begin{pf}
The idea is to apply the criterion of \cite[Proposition 2.4]{B}. Note that since $\beta$ and $\omega$ are rational, the image of $Z\colon \N(X)\to\C$ is a discrete subgroup of $\C$. Suppose one has a chain of monomorphisms in $\A=\A(\beta,\omega)$
\[\cdots\subset E_{i+1}\subset E_i\subset\cdots\subset E_1\subset E_0=E\]
with increasing phases $\phi(E_{i+1})>\phi(E_i)$ for all $i$. There are short exact sequences
\[0\lra E_{i+1}\lra E_i\lra F_i\lra 0\]
in $\A$. Since $\Im Z(A)$ is non-negative for all objects $A\in\A$
one has
\[\Im Z(E_{i+1})<\Im Z(E_i) \text{ for all }i.\]
Since the image of $Z$ is discrete it follows that for large enough $i$ the value of $\Im Z(E_i)$ is constant.
But then $\Im Z(F_i)=0$ so $\Re Z(F_i)\leq 0$ and so $\Re Z(E_{i+1})\geq \Re Z(E_i)$. But this contradicts $\phi(E_{i+1})>\phi(E_i)$.
Thus no such chain can exist.

Now suppose one has a chain  of epimorphisms
\[E=E_0\onto E_1\onto \cdots \onto E_{i}\onto E_{i+1}\onto\cdots\]
with decreasing phases $\phi(E_{i})>\phi(E_{i+1})$ for all $i$. There are short exact sequences
\[0\lra K_i\lra E_i\lra E_{i+1}\lra 0.\]
Applying the same argument as for the first part, and omitting a finite number of terms,
one can assume that $\Im Z(E_i)=\Im Z(E_{i+1})$ for all $i$. 
This time however this does not contradict the inequality
on the phases.

Taking long exact sequences in cohomology sheaves shows that there are epimorphisms
of sheaves
\[H^0(E_0)\onto H^0(E_1)\onto \cdots \onto H^0(E_{i})\onto H^0(E_{i+1})\onto\cdots.\]
Since the category of coherent sheaves is Noetherian this chain must terminate, and so, again omitting a finite number of terms, one can assume that
the map $H^0(E)\onto H^0(E_i)$ is an isomorphism for all $i$.

Consider the composite maps $E=E_0\onto E_i$ fitting into short exact sequences
\[0\lra L_i \lra E\lra E_{i}\lra 0.\]
Then there is a chain
\[0=L_0\subset L_1\subset\cdots\subset L_{i}\subset\cdots\subset E\]
and each $\Im Z(L_i)=0$ for all $i$. It will be enough to show that for large enough $i$ one has $L_i=L_{i+1}$.

There are short exact sequences
\[0\lra L_{i-1}\lra L_i\lra B_i\lra 0.\]
Taking cohomology sheaves shows that there are monomorphisms of sheaves
\[0=H^{-1}(L_0)\subset H^{-1}(L_1)\subset \cdots \subset H^{-1}(L_i)\subset \cdots\subset H^{-1}(E).\]
Once again, this chain must terminate,
so omitting a finite number of terms  one can assume
that the inclusion  $H^{-1}(L_i)\to H^{-1}(L_{i+1})$ is an isomorphism for all $i$.

For each $i$ there is a short exact sequence in $\A$
\[0\lra H^{-1}(L_i)[1]\lra L_i\lra H^0(L_i)\lra 0.\]
Since $\Im Z(L_i)=0$ one must also have $\Im Z(H^0(L_i))=0$. But by definition of the category $\A$
this is only possible if $H^0(L_i)$ is a torsion sheaf supported in dimension zero.
Taking cohomology gives a long exact sequence of sheaves
\[0\lra H^{-1}(B_i)\lra H^0(L_{i-1})\lra H^0(L_i)\lra H^0(B_i)\lra 0.\]
Since $H^{-1}(A)$ is torsion-free for any object $A\in\A$, and $H^0(L_{i-1})$ is torsion, one has $H^{-1}(B_i)=0$.

Now it will be enough to show that $H^0(B_i)=0$ for large $i$.
Equivalently one must bound the length of the finite-length sheaves $H^0(L_i)$.
Consider again the short exact sequence
\[0\lra L_i \lra E\lra E_{i}\lra 0.\]
Taking cohomology sheaves and recalling the assumption that the map $H^0(E)\to H^0(E_i)$ is an isomorphism for all $i$ gives long exact sequences
of sheaves
\[0\lra H^{-1}(L_i)\lra H^{-1}(E)\lRa{f} H^{-1}(E_i) \lra H^0(L_i) \lra 0.\]
Let $Q$ be the image of the map $f$. This is independent of $i$ up to isomorphism since the map $H^{-1}(L_i)\to H^{-1}(L_{i+1})$ is an isomorphism for all $i$.
But now there is a short exact sequence of sheaves
\[0\lra Q\lra H^{-1}(E_i)\lra H^0(L_i)\lra 0\]
in which the middle term $H^{-1}(E_i)$ is torsion-free by definition of the category $\A$. It follows that $Q$ is torsion-free and $H^0(L_i)$ is a subsheaf of the finite-length sheaf $Q^{**}/Q$. In  particular the length of $H^0(L_i)$ is bounded. This shows that $B_i=0$ for large $i$, and hence $E_i=E_{i+1}$ for large $i$, which contradicts the assumption of decreasing phase.

 Finally, the local-finiteness condition follows from Lemma \ref{disc} because the
image of the homomorphism $Z$ is a discrete subgroup of $\C$.
\end{pf}

% ***************************************************************************
% ***************************************************************************
% ***************************************************************************
% ***************************************************************************

\section{The covering map property}

Recall from the introduction that the subset $\Pol(X)\subset\N(X)\tensor\C$ is defined to be the set of vectors whose real and imaginary parts span positive definite
two-planes in $\N(X)\tensor\R$. Recall also that
\[\Delta(X)=\{\delta\in\N(X):(\delta,\delta)=-2\},\]
and that for each
$\delta\in\Delta(X)$ there is a corresponding complex hyperplane
\[\delta^{\perp}=\{\mho\in\N(X)\tensor\C:(\mho,\delta)=0\}\subset
\N(X)\tensor\C.\]
The hyperplane complement
\[\Pol_0(X)=\Pol(X)\setminus\bigcup_{\delta\in\Delta(X)}\delta^{\perp}.\]
is an open subset of $\Pol(X)$ (see Proposition \ref{cover} below).
The aim of this section is to use the deformation result of \cite{B}
to show that the map
\[\Pi\colon \Stab(X)\lra \N(X)\tensor\C\]
is a covering map over this subset.

Recall that a continuous map of topological spaces $f\colon S\to T$ is a covering map if every $t\in T$ has an open neighbourhood $t\in V\subset T$
such that the restriction of $f$ to each connected component of $f^{-1}(V)\subset S$ is a homoemorphism onto $V$. It is an immediate consequence
that if $S$ is non-empty and $T$ is connected then $f$ is surjective. 

\begin{lemma}
\label{ineq}
Let $\operatorname{\|\cdot\|}$ be a norm on
the finite dimensional vector space $\N(X)\tensor\C$ and let $(-,-)$ denote
the Mukai bilinear form on $\N(X)\tensor\C$. Take a vector
$\mho\in\Pol(X)$. Then there is a constant $r>0$ (depending on $\mho$) such that
\[|( u,v )|\leq r\,\| u\|\, |(\mho,v)|\]
for all $u\in \N(X)\tensor\C$ and all
$v\in \N(X)\tensor\R$ with $(v,v)\geq 0$. If moreover $\mho\in\Pol_0(X)$
then  one can choose $r>0$ so that the same inequality holds for
all $v\in \Delta(X)$.
\end{lemma}

\begin{pf}
The assumption $\mho\in\Pol(X)$ means that there is a basis
$e_1,\cdots,e_n$ of the real vector space
$\N(X)\tensor\R$ which is orthogonal with respect to the
Mukai inner product with
signature $(+1,+1,-1,\cdots,-1)$, and such that the real and imaginary parts of $\mho$ are
a basis for the subspace spanned by $e_1$ and $e_2$. Since all norms on a
finite-dimensional vector space are equivalent, one might as well take
\[\|v\|^2=v_1^2 +v_2^2 +\cdots +v_n^2,\]
where $v_i$ denotes the $i$-th component of a vector
$v\in\N(X)\tensor\R$ with respect to the basis $e_1,\cdots,e_n$.
Note that for vectors $u,v\in\N(X)\tensor\R$, the Cauchy-Schwartz
inequality gives
\[|(u,v)|=|u_1 v_1+u_2 v_2 -u_3 v_3 -\cdots- u_n v_n|\leq
\|u\|\|v\|.\]
Furthermore, since the real and imaginary parts of $\mho$ are a basis for the subspace spanned by $e_1$ and $e_2$,
there is a $k>0$ such that
for all $v\in\N(X)\tensor\R$ one has
\[ v_1^2+v_2^2\leq k |(\mho,v)|^2.\]
If a vector $v\in\N(X)\tensor \R$ satisfies $v^2\geq 0$ then
$\|v\|^2\leq2(v_1^2+v_2^2)$ so the result follows.
In the case when $v\in\Delta(X)$ one  has
\[\|v\|^2=2(v_1^2+v_2^2+1),\]
so it will be enough to check that
\[1\leq m|(\mho,v)|^2\]
for some $m>0$. In other words, 
as $v$ varies in $\Delta(X)$, the quantity
$|(\mho,v)|$ is bounded below by a positive real number.
If $|(\mho,v)|\leq1$ then $\|v\|^2\leq 2(k+1)$  and $v$ lies in a bounded subset of $\N(X)\tensor\R$.
But there are
only finitely many integral points in any bounded subset of
$\N(X)\tensor\R$, and the assumption 
$\mho\in\Pol_0(X)$ implies that all the $|(\mho,v)|$ are nonzero, so the
result follows.
\end{pf}

It is worth recording the following consequence of the above proof.

\begin{lemma}
\label{em}
If $\mho\in\Pol_0(X)$ and $m>0$ is a constant then there are only finitely many elements $v\in \N(X)$ such that
$v^2\geq -2$ and $|(\mho,v)|\leq m$.
\end{lemma}

\begin{pf}
The inequalities given in the proof of Lemma \ref{ineq} give an upper bound for $\|v\|$, and there are only finitely many integral points in any bounded subset of $\N(X)\tensor\R$.
\end{pf}

The following is the main result of this section.

\begin{prop}
\label{cover}
The subset $\Pol_0(X)\subset\N(X)\tensor\C$ is open
and the restriction
\[\Pi^{-1}\,(\Pol_0(X))\lRa{\Pi} \Pol_0(X)\]
is a covering map.
\end{prop}

\begin{pf}
Fix a norm $\|\cdot\|$ on
the finite-dimensional complex vector space
$\N(X)\tensor\C$ and take a point $\mho\in\Pol_0(X)$.
Lemma \ref{ineq} below shows that there is a constant $r>0$ (depending on $\mho$)
such that
\[|( u,v )|\leq r\,\| u\|\, |(\mho,v)|\]
for all $u\in \N(X)\tensor\C$ and all
$v\in \N(X)\tensor\R$ with $v^2\geq 0$ and all
$v\in \Delta(X)$.
Given $\eta>0$, define an open subset
\[B_{\eta}(\mho)=\{\mho'\in\N(X)\tensor\C:\|\mho'-\mho\|
<\eta/r\}\subset\N(X)\tensor\C.\]
Then one has an implication
\[\mho'\in B_{\eta}(\mho)\implies |(\mho',v)-(\mho,v)|<\eta\,|(\mho,v)|\]
for all $v\in\N(X)\tensor\R$ with $v^2\geq 0$ and all
$v\in\Delta(X)$. It follows from this that if $\eta<1$ then
any $\mho'\in B_{\eta}(\mho)$ spans a positive definite two-plane in $\N(X)\tensor\R$
and satisfies $(\mho',v)\neq 0$ for all $v\in \Delta(X)$.
Hence $B_{\eta}(\mho)\subset\Pol_0(X)$ and, in particular, $\Pol_0(X)$ is an open
subset of $\N(X)\tensor\C$.

For each $\sigma\in\Stab(X)$ with $\Pi(\sigma)=\mho$
define an open subset
\[C_{\eta}(\sigma)=\{\tau\in\Pi^{-1}(B_{\eta}(\mho)):
f(\sigma,\tau)<\ha\}\subset\Stab(X),\]
where $f$ is the function introduced in Section \ref{two}.
Lemma \ref{cleave} shows that for any
$\mho'\in B_{\eta}(\mho)$ and any stable object
$E\in\D(X)$
\[|(\mho',v(E))-(\mho,v(E))|<\eta\,|(\mho,v(E))|.\]
Thus by Theorem \ref{biggy},
for small enough $\eta>0$ the map
\[\Pi\colon C_{\eta}(\sigma)\lra B_{\eta}(\mho)\]
is onto, and hence, by Theorem \ref{cath} and Lemma
\ref{rachael}, a homeomorphism. It follows from this that
every stability condition in $\Pi^{-1}(\Pol_0(X))$ is full (see Definition \ref{full}).

Fix a positive real number $\epsilon<\ei$ and
assume that $\eta<\ha\sin(\pi\epsilon)$. Then, by Theorem \ref{biggy} again, and using Lemma \ref{fully},
for each $\sigma\in\Pi^{-1}(\mho)$ the subset $C_{\eta}(\sigma)$ is mapped
homeomorphically by $\Pi$ onto $B_{\eta}(\mho)$. The final thing to
check is that one has a disjoint union
\[\Pi^{-1}(B_{\eta}(\mho))=\bigcup_{\sigma\in\Pi^{-1}(\mho)}C_{\eta}(\sigma).\]

The fact that the union is disjoint follows from Lemma \ref{rachael}.
Suppose $\tau$ is a stability condition on $\D(X)$
with $\Pi(\tau)=\mho'\in B_{\eta}(\mho)$. If an object $E\in\D(X)$
is stable in $\tau$
then
\[|(\mho',v(E))-(\mho,v(E))|<\eta\,|(\mho,v(E))|
<\frac{\eta}{1-\eta}\,|(\mho',v(E))|<2\eta\,|(\mho',v(E))|.\]
Applying Theorem \ref{biggy} again
gives a stability condition $\sigma\in\Pi^{-1}(\mho)$ such
that $f(\sigma,\tau)<\epsilon$, and then $\tau\in C_{\eta}(\sigma)$.
\end{pf}

The following definition will be useful.

\begin{defn}
\label{han}
A connected component $\Stab^*(X)\subset\Stab(X)$ will be called good if it contains a point $\sigma$ with
$\pi(\sigma)\in \Pol_0(X)$. A stability condition $\sigma\in\Stab(X)$ will be called good if it lies in a good component.
\end{defn}

As observed in the proof of Proposition \ref{cover}, a good component of $\Stab(X)$ is also full. Indeed, if a connected component $\Stab^*(X)\subset\Stab(X)$ is good then by Proposition \ref{cover}
the image of the map
\[\pi\colon \Stab^*(X)\lra \N(X)\tensor\C\]
contains one of the two connected components of the open subset $\Pol_0(X)$. In particular, this image is not contained in a linear subspace.

The results of Sections 6 and 7 show that there do indeed exist good components of $\Stab(X)$.
Note also that if $\Phi\in\Aut\D(X)$
is an autoequivalence then it acts on $\N(X)$ as an isometry, so if
a stability condition $\sigma\in\Stab(X)$ satisfies $\pi(\sigma)\in\Pol_0(X)$ then the same is
true of $\Phi(\sigma)$. Thus $\Aut\D(X)$ acts on the set of good stability conditions.

% ***************************************************************************
% ***************************************************************************
% ***************************************************************************
% ***************************************************************************

\section{The wall and chamber structure}
\label{wall}

Consider for a moment
 the problem of determining the set of $\mu_{\omega}$-stable
sheaves on $X$ as a function of the ample divisor class
$\omega\in\Amp(X)$. As is well-known, if one restricts attention to sheaves
with fixed Mukai vector, the space $\Amp(X)$
splits into a series of chambers in such a way that the set of
$\mu_{\omega}$-stable sheaves with the given Mukai vector is constant in each
chamber.
The aim of this section is to prove a similar result for $\Stab(X)$.

Suppose $\Stab^*(X)\subset\Stab(X)$ is a good component. The main result of this section is that
if $\S$ is some
finite set of objects of $\D(X)$, and $B\subset\Stab^*(X)$
is a compact subset, then there
is a finite collection of walls,
which is to say
codimension
one submanifolds of $B$, such that as one varies the stability
condition $\sigma\in B$, an element of $\S$
which is stable
can only become unstable if one crosses a wall.

More generally it is not necessary to assume that the set
$\S\subset\D(X)$ is finite, only
that its elements have bounded mass in the sense of the following definition.

\begin{defn}
A set of objects $\S\subset\D(X)$ has bounded
mass in a connected component $\Stab^*(X)\subset\Stab(X)$ if
\[\sup\{\m_{\sigma}(E): E\in \S\}<\infty\]
for some point $\sigma\in \Stab^*(X)$.
Note that the fact that $\Stab^*(X)\subset\Stab(X)$ is connected implies that
$d(\sigma,\tau)<\infty$ for all points $\sigma,\tau\in\Stab^*(X)$, so that if
this condition holds at some point $\sigma\in\Stab^*(X)$ then it holds at
all points.
\end{defn}

The following easy result will be crucial.

\begin{lemma}
Suppose the subset $\S\subset\D(X)$ has bounded mass in a good component $\Stab^*(X)\subset \Stab(X)$.
Then the set of Mukai vectors $\{v(E):E\in \S\}$ is finite.
\end{lemma}

\begin{pf}
By assumption there is a $\sigma=(Z,\P)\in\Stab^*(X)$ such that
$\Pi(\sigma)\in\Pol^+_0(X)$.
Let $m>0$ be such that $m_{\sigma}(E)<m$ for all $E\in \S$.
Then the stable factors $A_1,\cdots,A_n$ of an object $E\in\S$ with
respect to the stability condition $\sigma$ must satisfy
\[\sum_i|Z(A_i)|<m.\]
By Lemma \ref{cleave} and Lemma \ref{em} there are only finitely many 
possibilities for the Mukai vectors $v(A_i)$, and hence
only finitely many possibilities for the Mukai vector $v(E)$.
\end{pf}

The next result gives the claimed wall and chamber structure.

\begin{prop}
\label{lo}
Suppose the subset $\S\subset\D(X)$ has bounded mass in a good component $\Stab^*(X)\subset\Stab(X)$,
 and
fix a compact subset $\co\subset\Stab^*(X)$.
Then there is a
finite collection $\{W_{\gamma}:\gamma\in\Gamma \}$
of (not necessarily closed) real codimension one submanifolds of
$\Stab^*(X)$ such that any connected component
\[\CC\subset \co\setminus\bigcup_{\gamma\in\Gamma}W_\gamma\]
has the following property:
if $E\in \S$ is semistable in $\sigma$
for some $\sigma\in\CC$, then
$E$ is semistable in $\sigma$ for all $\sigma\in\CC$;
if moreover $E\in \S$ has primitive Mukai vector then $E$ is 
stable in $\sigma$ for all $\sigma\in\CC$.
\end{prop}

\begin{pf}
Let $\TT\subset\D(X)$ be the set of nonzero objects $A\in \D(X)$
such that for some
$\sigma\in\co$ and some $E\in\S$ one has $m_{\sigma}(A)\leq m_{\sigma}(E)$.
The fact that $\co$ is compact
implies that the
quotient $m_{\tau}(E)/m_{\sigma}(E)$ is uniformly bounded for all
nonzero $E\in\D(X)$, and for all $\sigma,\tau\in \co$, so
the subset $\TT\subset\D(X)$ has bounded mass in $\Stab^*(X)$.
Note that if $A$ is a semistable factor of
an object $E\in\S$ in some stability condition $\sigma\in\co$ then $A$
is an element
of $\TT$.

Let $\{v_i:i\in I\}$ be the finite set of Mukai vectors of objects of
$\TT$ and let $\Gamma$ be the set of pairs $i,j\in I$ such that $v_i$
and $v_j$ do not lie on the same real line in $\N(X)\tensor\R$.
For each $\gamma\in\Gamma$ define
\[W_{\gamma}=\{\sigma=(Z,\P)\in\Stab^*(X):Z(v_i)/Z(v_j)\in
\R_{>0}\}.\]
Since $\Stab^*(X)$ is a good component, and hence full,
the map
$\Pi\colon\Stab^*(X)\to\N(X)\tensor\C$
is a local homeomorphism. Since $W_{\gamma}$ is
the inverse image under $\Pi$ of an open
subset of a real quadric in
$\N(X)\tensor\C$, it follows that each
$W_{\gamma}$ is a real codimension one submanifold of $\Stab^*(X)$.

Suppose that $\CC\subset \co$ is a connected component of
\[\co\setminus\bigcup_{\gamma\in\Gamma}W_\gamma\]
and fix an object $E\in\S$. Consider the subset $V\subset\CC$
consisting of points at which $E$ is semistable, and assume that $V$
is nonempty.
The definition of the topology on $\Stab(X)$ ensures that $V$ is a closed
subset of $\CC$.
The next step is to show that $V$ is also open in $\CC$ and hence that
$V=\CC$.

Suppose then that $\sigma=(Z,\P)\in V$ is such that $E\in\P(\phi)$
for some $\phi\in\R$. Take $0<\eta<\ei$ such that the open neighbourhood
\[U=\{\tau\in\Stab(X):d(\sigma,\tau)<\eta\}\]
is contained in $\CC$. It will be enough to show that $U\subset V$.

Note that the assumption $\eta<\ei$ ensures that if  $A$ is a semistable factor of $E$ in some stability condition in $U$ then $A$ lies in the abelian subcategory
\[\A=\P\big((\phi-\ha,\phi+\ha]\big)\subset\D(X),\]
and that moreover, for any other stability condition $\tau=(W,\Q)\in U$ the central charge $W(A)$ lies in the half-plane
\[H_{\phi}=\big\{r\exp(i\pi\psi):r>0\text{ and } \phi-\ha<\psi<\phi-\ha\big\}.\] 
Suppose for a contradiction that
$E$ is unstable at some point
$\sigma'=(Z',\P')\in U$.
Then there is a semistable factor $A$ of $E$ in the
stability condition $\sigma'$ which
is a subobject of $E$ in the category $\A$,
and which satisfies
$\Im Z'(A)/Z'(E)>0$. As $\tau=(W,\Q)$ varies in $U$ the complex
numbers $W(A)$ and $W(E)$ remain in $H_{\phi}$, and
therefore, since $U$
is contained in a connected component of \[\co\setminus\bigcup_{\gamma\in\Gamma}W_\gamma\]
and $A$ and $E$ are objects of $\TT$, it follows that
$\Im W(A)/W(E)>0$ for all $\tau=(W,\Q)\in U$. This contradicts the fact that $E$ is semistable at
$\sigma$ thus proving the claim.

To complete the proof of the Proposition,
suppose an object $E\in\S$ has primitve Mukai vector
$v(E)\in\N(X)$, and suppose
$E$ is semistable but not stable at some point
$\sigma\in \co$. Each stable factor $A$ of $E$ in the stability
condition $\sigma$ has mass less than $E$ so has Mukai vector $v_i$
for some $i\in I$.
Since $v(E)$ is primitive, not all the Mukai vectors of the stable
factors of $E$ are multiples of the Mukai vector of $E$. But the phases of
all the stable factors are the same so
one must have $\sigma\in W_{\gamma}$ for some pair
$\gamma=(i,j)\in\Gamma$.
\end{pf}

\begin{prop}
\label{open}
Suppose the subset $\S\subset\D(X)$ has bounded mass in a good component $\Stab^*(X)\subset\Stab(X)$ and that
each object $E\in\S$ has primitive Mukai vector $v(E)\in\N(X)$.
Then the set of points $\sigma\in\Stab^*(X)$ at
which all objects of $\S$ are stable is open in $\Stab^*(X)$.
\end{prop}

\begin{pf}
It will be enough to fix a compact subset
$\co\subset\Stab^*(X)$
and prove that the set
$F=\{\sigma\in\co:\text{not every }E\in\S
\text{ is stable in }\sigma\}$
is a closed subset of $\co$.
Let $\TT$ by the set of all semistable factors of objects of
$\S$ in all stability conditions of $\co$. It was observed in the proof of Proposition \ref{lo} that
this set has bounded mass in $\Stab^*(X)$.
Consider
the corresponding chamber decomposition given by Proposition \ref{lo}.
The aim is to show that $F$ is closed by proving
that it is the union of the closures of those chambers $\CC$
in which some object $E\in\S$ is not stable.

Take a stability condition $\sigma=(Z,\P)\in F$ lying in the closure
of a finite set of chambers $\CC_j\subset\co$. Take an object
$E\in\S$ such that $E$ is not
stable in $\sigma$.
If $E$ is
not semistable in $\sigma$ then $E$ is not semistable in an
open neighbourhood of
$\sigma$, so $E$ is not stable in each chamber $\CC_j$.
The other possibility is that $E$ is
semistable of some phase $\phi$ in $\sigma$
but not stable.
Then there is a short exact sequence
\[0\lra A\lra E\lra
B\lra 0\]
in the abelian subcategory $\P(\phi)\subset\D(X)$.
Since $v(E)$ is primitive, the Mukai vectors of $A$
and $B$ are not multiples of each other, and since $\Pi$ is a local
homeomorphism it follows that
there are points $\tau=(W,Q)$ arbitrarily close to $\sigma$ for which $\Im W(A)/W(B) >0$. This implies that
$E$ must be unstable in one of the chambers $\CC_j$.
Thus either way $\sigma$ lies in the closure of a
chamber $\CC\subset\co$ in which $E$ is unstable.

For the converse suppose that some object $E\in\S$ is unstable in a chamber
$\CC\subset\co$. Take a stability condition
$\sigma\in\CC$ and let the semistable factors of $E$ in
$\sigma$ be $A_1,\cdots ,A_n$ with phases
$\phi(A_1)>\cdots>\phi(A_n)$.
Then each $A_i$ is an object of $\TT$, so by the construction of
Proposition \ref{lo},
the objects $A_i$ are
semistable at each point of the closure of $\CC$ and
satisfy $\phi(A_1)\geq\cdots\geq\phi(A_{n})$. It follows that $E$ is
not stable at any point of the closure of $\CC$, which is to say that
$F$ contains the closure of $\CC$.
\end{pf}

% ***************************************************************************
% ***************************************************************************
% ***************************************************************************
% ***************************************************************************

\section{Classifying stability conditions}
\label{class}

It follows from Lemma \ref{balloon} that all the stability conditions
constructed in Section \ref{con} have the property that for
each point $x\in X$, the
corresponding skyscraper sheaf $\OO_x$ is stable in $\sigma$ of phase one. The next step is to
prove a  converse to this result.

\begin{lemma}
\label{lemma}
Suppose $\sigma=(Z,\P)\in\Stab(X)$ is a stability condition on $X$
such that for each point $x\in X$ the sheaf
$\OO_x$ is stable in $\sigma$ of phase one. Let $E$ be an object
of $\D(X)$. Then
\begin{itemize}
\item[(a)] if $E\in\P((0,1])$
then the cohomology sheaves $H^i(E)$ vanish unless
$i\in\{-1,0\}$, and moreover the sheaf $H^{-1}(E)$ is torsion-free,

\item[(b)] if $E\in\P(1)$ is stable then either $E=\OO_x$ for some
$x\in X$, or $E[-1]$ is a locally-free sheaf,

\item[(c)] if $E\in\Coh(X)$ is a sheaf  then $E\in\P((-1,1])$; if $E$ is a torsion sheaf then $E\in\P((0,1])$,

\item[(d)] the pair of subcategories
\[\mathcal{T}=\Coh(X)\cap \P((0,1])\text{ and }
\mathcal{F}=\Coh(X)\cap \P((-1,0])\]
defines a torsion pair on the category of coherent sheaves $\Coh(X)$
and moreover the category $\P((0,1])$
is the corresponding tilt.
\end{itemize}
\end{lemma}

\begin{pf}
Suppose that $E$ is stable of phase
$\phi$ for some $0< \phi< 1$.
For any $x\in X$
the sheaf $\OO_x$ is stable of phase 1, so that $\Hom^i_X(E,\OO_x)=0$
for $i<0$, and Serre duality gives
\[\Hom^i_X(E,\OO_x)=\Hom^{2-i}_X(\OO_x,E)^{\dual}=0\]
for $i\geq 2$.
Taking a locally-free resolution of $E$ and truncating (see for example \cite[Proposition 5.4]{B2})
it follows that $E$ is isomorphic to a length two complex of locally-free
sheaves, 
so that $E$ satisfies the conclusions of part (a).

If $E$ is stable of phase $\phi=1$ and $E$ is not a skyscraper sheaf, then there cannot be any nonzero maps $E\to \OO_x$ or $\OO_x \to E$, so
the same argument shows that
$\Hom^i_X(E,\OO_x)=0$ unless $i=1$, and hence $E[-1]$ is
locally-free.
Any object of $\P((0,1])$ has a filtration by stable objects with
phases in the interval $(0,1]$ so this proves (a) and (b).

Suppose now that $E$ is a sheaf on $X$.
For any object $A\in\P(\smash{>}{1})$ part (a) shows that
$H^i(A)=0$ for $i\geq 0$ so
that $\Hom_X(A,E)=0$. Similarly, if $B\in\P(\smash{\leq}{-1})$ then
$H^i(B)=0$ for $i\leq 0$ so that $\Hom_X(E,B)=0$. It follows that
$E\in\P((-1,1])$ which gives the first half of (c).

Given a sheaf $E$ on $X$, by the first part of (c) there is a triangle
\begin{equation*}
\xymatrix@C=.5em{ D \ar[rrrr] &&&& E  \ar[dll] \\
 && F \ar@{-->}[ull]} \end{equation*}
with $D\in\P((0,1])$ and $F\in\P((-1,0])$. By part (a)
one has $H^i(D)=0$ unless $i\in\{-1,0\}$ and $H^i(F)=0$ unless
$i\in\{0,1\}$. Taking the long exact sequence in cohomology
one concludes that $D$ and $F$ must both be sheaves.
This shows that $\mathcal{T},\mathcal{F}$ is a torsion pair and it is immediate that $\P((0,1])$ is the corresponding tilt.
Moreover, by part (a) again
the sheaf $F$ is torsion-free so that if $E$ is torsion then
$E\in\P((0,1])$, which completes the proof of (c).
\end{pf}

Certain technical problems arise if one tries to classify all stability conditions satisfying the conditions of Lemma \ref{lemma}.
However the extra assumption that $\sigma$ is good (Definition \ref{han}) leads to a complete classification.

\begin{defn}
Let $U(X)\subset\Stab(X)$ be the subset consisting of  good stability conditions $\sigma\in\Stab(X)$
such that for each point $x\in X$ the sheaf
$\OO_x$ is stable in $\sigma$ of the same phase.
\end{defn}

It follows from Proposition \ref{open} that $U(X)\subset\Stab(X)$ is open.

Recall from the introduction that $\GL(2,\R)$ acts freely on $\Pol(X)$ by change of basis
in the positive two plane spanned by the real and imaginary parts of a vector $\mho\in\N(X)\tensor\C$.
A section of this action is provided by
the submanifold
\[\Q(X)=\big\{\mho\in \Pol(X):(\mho,\mho)=0,\
(\mho,\bar{\mho})>0,\  r(\mho)=1\big\}\subset\N(X)\tensor\C,\]
where $r(\mho)$ is the projection of $\mho\in\N(X)\tensor\C\subset H^*(X,\C)$ onto $H^0(X,\C)$.
Note that $\Q(X)$ can be identified with the tube domain
\[\big\{\beta+i\omega\in \NS(X)\tensor\C:
\omega^2>0\big\}\]
via the exponential map
$\mho=\exp(\beta+i\omega)$.

Before stating the next result, let us agree to say that a stability condition $\sigma=(Z,\P)$ on $\D(X)$ arises from
the construction of Section \ref{con} if there is a pair $\beta,\omega\in\NS(X)\tensor\R$ with $\omega\in\Amp(X)$ such that
the central charge $Z$ is given by
\[Z(E)=(\exp(\beta+i\omega),v(E)),\]
and that moreover the
heart $\P((0,1])$ of $\sigma$ is the abelian subcategory $\A(\beta,\omega)\subset\D(X)$ defined in Section \ref{con}.
Note that in this case the pair $\beta,\omega$ is uniquely defined by $\sigma$, and that conversely,
by Proposition \ref{pg}, given any such pair there is at most one stability condition arising from $\beta,\omega$ via the construction of Section \ref{con}.

\begin{prop}
\label{tilt}
If $\sigma\in U(X)$ 
then there is a unique element $g\in\grp$,
such that $\sigma g$ arises from the construction of Section \ref{con}.
\end{prop}

\begin{pf}
The proof is best broken into a sequence of steps.

\Step1
Applying an element of $\grp$ one can assume that $\OO_x\in\P(1)$ for all $x\in X$.
Since $Z(\OO_x)$ lies on the real axis, the central charge $Z$
of an object $E\in\D(X)$ with Mukai vector $(r,\Delta,s)$ must
satisfy
\[\Im Z(E)=\Delta\cdot\omega -rx, \]
for some $\omega\in\NS(X)\tensor\R$ and some $x\in\R$.
The first step is to show that $\omega\in\Amp(X)$ is ample.

Suppose $C\subset X$ is a curve.
Lemma \ref{lemma}(c) shows that the torsion sheaf  $\OO_C$ lies in the subcategory $\P((0,1])$.
If $Z(\OO_C)$ lies on the real axis it follows
that $\OO_C\in\P(1)$ which is impossible by Lemma \ref{lemma}(b).
Thus $\Im Z(\OO_C)=C\cdot \omega >0$.
Now $\sigma$ is a good stability condition, and hence full,  so the map
\[\pi\colon\Stab(X)\to\N(X)\tensor\C\]
is a local homeomorphism near $\sigma$. Since $U(X)$ is open
 one may deform $\omega$ a little and still conclude that $C\cdot \omega>0$ for all curves $C\subset X$. Thus
 $\omega$ lies in the interior of the nef cone in $H^2(X,\R)$ and hence is ample.

\Step2
Consider the torsion pair $(\mathcal{T},\mathcal{F})$ of Lemma \ref{lemma}(d).
Note that by Lemma \ref{lemma}(c) all torsion sheaves lie in $\mathcal{T}$.
The next claim is that 
if $E$ is a $\mu_\omega$-stable torsion-free sheaf then either $E\in\mathcal{T}$ or $E\in\mathcal{F}$ depending
on whether $\Im Z(E)> 0$ or $\Im Z(E)\leq 0$.

Given such an $E$ there is a short
exact sequence of sheaves
\[0\lra D\lra E\lra F\lra 0\]
with $D\in\mathcal{T}$ and $F\in\mathcal{F}$. Assume $D$ and $F$ are both nonzero.
By Lemma \ref{lemma}(c), the sheaf
$F$ must be torsion-free. Now one has
$\mu_\omega(D)\geq x$ and
$\mu_\omega(F)\leq x$ and since $E$ is
$\mu_\omega$-stable, this gives a contradiction. Thus $E$ lies in either $\T$ or
$\F$. Clearly,
if $\Im Z(E)>0$ then $E\in\T$ and if
$\Im Z(E)<0$ then $E\in\F$.

Suppose that $\Im Z(E)=0$. It must be shown that $E\in\F$. 
Suppose for a contradiction that $E\in\T$. Since $Z(E)\in\R$
it follows that $E\in\P(1)$.
Consider a nonzero
map $f\colon E\to\OO_x$ and let $F$ be its kernel in the category
$\Coh(X)$. Since $\OO_x$ is
stable, the map $f$ is a
surjection in the abelian category $\P(1)$, so the object $F$ also lies in
$\P(1)$ and hence in $\T$.
Repeating must eventually give $Z(F)\in\R_{>0}$ which contradicts
$F\in\P((0,1])$.

\Step3
The next step is to show that $\mho=\pi(\sigma)\in\Pol(X)$, which is
to say that the real and imaginary parts of $\mho$ span a positive definite two-plane in $\N(X)\tensor\R$.
Equivalently one must show that $\Ker(Z)\subset \N(X)\tensor \R$ is negative definite.

Suppose for a contradiction that there is a $v\in\NS(X)\tensor\R$ such that $Z(v)=0$ and $(v,v)\geq 0$.
Since $U(X)$ is open one can deform $\sigma$ and assume that $(v,v)>0$ and $v\in\NS(X)\tensor\QQ$,
and hence by removing denominators, that $v\in\NS(X)$.

If $r(v)=0$ then the assumption $(v,v)>0$ implies that either $v$ or $-v$ is the Mukai vector of a torsion sheaf $T$ on $X$. But by Lemma \ref{lemma}(c) one  has $T\in\P((0,1])$ so that $Z(v)\neq 0$. Otherwise one can assume that $r(v)>0$. But Theorem \ref{yosh} then shows that $v$ is the Mukai vector of a $\mu_{\omega}$-semistable torsion-free sheaf $E$. By Step 2 it follows that either $E\in\mathcal{T}$ or $E\in\mathcal{F}$ and hence  $Z(v)\neq 0$.

\Step4
Now one may apply an element of $g\in \grp$ and assume that $\mho\in \Q(X)$ so that
\[Z(E)=(\exp(\beta+i\omega),v(E)),\]
for some pair $\beta,\omega\in\NS(X)\tensor\R$
with $\omega^2>0$. The  group $\GL(2,\R)$ acts freely on $\Pol(X)$ and
the submanifold $\Q(X)$ is a section of this action, so the element $g$ is defined uniquely up to the kernel of the homomorphism $\grp\to\GL(2,\R)$. This kernel acts on $\Stab(X)$ by even shifts, so if one also assumes that each object $\OO_x$ is stable in $\sigma g$ of phase 1 then $g$ is uniquely determined.

By Step 1 one knows that $\omega$ is ample. Then Step 2 shows that the torsion pair $(\mathcal{T},\mathcal{F})$ of Lemma \ref{lemma}(d)
coincides with the torsion pair of Lemma \ref{mckay}.
It follows immediately that $\P((0,1])=\A(\beta,\omega)$ and hence that $\sigma$ arises from the construction of Section \ref{con}.
\end{pf}

% ***************************************************************************
% ***************************************************************************
% ***************************************************************************
% ***************************************************************************

\section{The open subset $\U$}
\label{U}

Following Proposition \ref{tilt} it is now possible to give a precise description of the set $U(X)\subset \Stab(X)$
defined in the last section.
Firstly one knows that the action of $\grp$ on $U(X)$ is free.
A section of this action is defined by the submanifold
\[\Ured=\big\{\sigma\in U(X):\Pi(\sigma)\in\Q(X)
\text{ and each }\OO_x \text{ is stable in $\sigma$ of
phase 1}\big\}.\]
The proof of Proposition \ref{tilt} showed that $\Ured$ is precisely the set of stability conditions arising from the construction of Section \ref{con}. In particular, the map $\pi$ is injective when restricted to $\Ured$.

To understand the image of  $\Ured$
recall first that $\Q(X)$ can be identified with the tube domain
\[\big\{\beta+i\omega\in \NS(X)\tensor\C:
\omega^2>0\big\}\]
via the exponential map
$\mho=\exp(\beta+i\omega)$.
The image of the complexified ample cone under this map
defines an open subset
\[\KK(X)=\big\{\mho=\exp(\beta+i\omega)\in\Q(X):\omega\in\Amp(X)\big\}\subset\Q(X).\]
Define  \[\Delta^+(X)=\big\{\delta\in\N(X):(\delta,\delta)=-2\text{ and }r(\delta)>0\big\}\subset\Delta(X).\]
Proposition \ref{below} below shows that the image of $\Ured$
under $\pi$ is the subset
\[\KKo(X)=\big\{\mho\in\KK(X):(\mho,\delta)\notin\R_{\leq 0}\text{ for
all }\delta\in\Delta^+(X)\big\}\subset \Q(X).\]
Note that this is in fact a subset of $\P_0(X)$ since if $\delta\in\Delta(X)$
satisfies $r(\delta)=0$ then $\delta=\pm(0,C,n)$ for some $(-2)$-curve $C\subset X$ and it follows that
$\Im (\mho,\delta)\neq 0$ for any $\mho\in\KK(X)$.

\begin{lemma}
The subset $\KKo(X)\subset \Q(X)$ is open and contractible.
\end{lemma}

\begin{pf}
To show that $\KKo(X)$ is open note that it is the complement in
$\KK(X)$ of the real
halfplanes
\[ H(\delta)=\big\{\mho\in\KK(X):(\mho,\delta)\in\R_{\leq 0}\big\}\]
for elements $\delta\in\Delta^+(X)$. Thus it will be enough to show that
these hyperplanes are locally finite in $\KK(X)$.

Fix a bounded region $\co\subset\KK(X)$ and take an element
$\mho=\exp(\beta+i\omega)\in \co$. Suppose that
$\mho\in H(\delta)$ for some
$\delta=(r,\Delta,s)\in\Delta^+(X)$.
According to the formula $(\star)$ of Section \ref{con}
\[\frac{1}{2r}\bigg((\Delta^2-2rs)+r^2\omega^2-(\Delta-r\beta)^2\bigg)
+i(\Delta-r\beta)\cdot\omega=(\mho,\delta)\in\R_{\leq 0}\]
Thus $\Delta-r\beta$ is an element of
the sublattice $\omega^{\perp}\subset\NS(X)\tensor\R$, which is
negative definite.
Now $\Delta^2-2rs=-2$, so that
\[r^2\omega^2-(\Delta-r\beta)^2\leq 2,\]
and since $\omega$ is constrained to lie
in a bounded region of $\NS(X)\tensor\R$,
it follows that there are only finitely many possible choices
for $r$. Also $\Delta-r\beta$ lies in a bounded region of $\omega^{\perp}\subset\NS(X)\tensor \R$, and since $\beta$ is constrained to lie in a bounded region of $\NS(X)\tensor \R$ it follows that there are only finitely many possibilities for $\Delta$, and hence
only finitely many hyperplanes $H(\delta)$ which meet
the region $\co$.

To see that $\KKo(X)$ is contractible note that
if $\exp(\beta+i\omega)$ is an element of $\KKo(X)$ then so
is $\exp(\beta+it\omega)$ for every $t>1$.
Thus $\KKo(X)$ has a deformation retraction onto
the subset
\[\big\{\exp(\beta+i\omega)\in\KK(X):\omega^2> 2\big\}\subset\KKo(X)\]
and hence is contractible.
\end{pf}

The following result is a stronger version of Proposition \ref{tilt}.

\begin{prop}
\label{below}
The map $\Pi$ restricts to
give a homeomorphism
\[\Pi\colon\Ured\lra\KKo(X).\]
The stability condition
in $\Ured$ corresponding to a point $\exp(\beta+i\omega)\in \KKo(X)$
is given by the construction of Section \ref{con}.
\end{prop}

\begin{pf}
As noted above, any stability condition $\sigma\in\Ured$ arises from the construction of Section \ref{con},
so the map $\Pi$ is injective when restricted
to  $\Ured\subset\Stab(X)$ and
hence gives a homeomorphism $\Ured\to\Pi(\Ured)$.
Moreover, for any stability condition $\sigma\in\Ured$
one has $\Pi(\sigma)=\exp(\beta+i\omega)\in\KK(X)$.
Theorem \ref{yosh} ensures that for
any $\delta\in\Delta^+(X)$ there is a $\mu_{\omega}$-semistable
torsion-free sheaf $E$ with $v(E)=\delta$.
Since $\sigma$ arises from the construction of Section \ref{con} one either has $E\in \P((0,1])$ or $E[1]\in \P((0,1])$. Either way, $Z(E)\neq 0$. This shows that $\pi(\sigma)\in\KKo(X)$.
Thus $\Pi$ maps $\Ured$ into $\KKo(X)$, and it only remains to show
that this map is onto.
But Proposition \ref{open} shows that $\Pi(\Ured)$ is an open subset of
$\KKo(X)$, so since $\KKo(X)$ is connected it will be enough to show that
it is also closed.

Suppose $\mho=\exp(\beta+i\omega)\in\KKo(X)$ lies in the
 closure of $\Pi(\Ured)$. Since $\Pi$ is a covering over $\KKo(X)$
 there is a stability condition $\sigma$ lying in the closure of
 $\Ured$ with $\Pi(\sigma)=\mho$. If $\sigma$ is not an element of
$\Ured$
there must be
some skyscraper sheaf $\OO_x\in\P(1)$ which is semistable in $\sigma$
 but not
stable.
Since $\omega$ is ample, any class $v\in\N(X)$ such that $\Im Z(v)=0$ must either  be a multiple
of $v(\OO_x)$ or have nonzero rank. Thus there must be a stable factor $A$ of $\OO_x$ with $r(A)>0$.
Then $Z(A)$ lies on the negative real axis. The argument of
Lemma \ref{slopefun} shows that $A$ must be spherical, so that $v(A)\in\Delta^+(X)$,
and since $\mho\in\KKo(X)$ this is impossible. Thus $\sigma\in\Ured$ and $\Pi(\Ured)$ is closed in $\KKo(X)$.
\end{pf}

Recall that $\Pol^+(X)\subset\Pol(X)$ is the connected component containing $\KK(X)$, and that $\Pol^+_0(X)=\Pol^+(X)\cap\Pol_0(X)$. It is worth recording here the following easy consequence of what has been proved so far.

\begin{cor}
\label{liz}
If $\sigma\in \U$ then $\pi(\sigma)\in\Pol^+_0(X)$. Moreover $\sigma$ is determined by $\pi(\sigma)$ up to even shifts.
\end{cor}

\begin{pf}
As shown in Proposition \ref{tilt}, for $\sigma\in\U$ there is a unique $g\in\grp$ such that $\sigma g\in\Ured$.
The first statement then follows from Proposition \ref{below} since $\Pol^+_0(X)$ is invariant under $\GL(2,\R)$.
For the second statement note that the only elements of the group $\grp$ which fix the central charge of a stability condition are double shifts.
\end{pf}

Note that Proposition \ref{below} 
shows that the stability
function of Lemma \ref{slopefun} always has the Harder-Narasimhan property.
Another consequence of Proposition \ref{below}
is that the set $\Ured$, and hence also $U(X)$, is connected. Thus all the stability
conditions constructed in Section \ref{con} are contained in the same
connected component of $\Stab(X)$.

\begin{defn}
Let $\Stab^\dagger(X)\subset\Stab(X)$ be the unique connected component
containing the set $U(X)$.
\end{defn}

The results of this section give a precise description
of the stability conditions lying in
the open subset $\U\subset\Stab(X)$. The next step is to analyse the
boundary of this set.

% ***************************************************************************
% ***************************************************************************
% ***************************************************************************
% ***************************************************************************

\section{The boundary of $\U$}
\label{dU}

The aim of this section is to study the boundary
$\partial U(X)=\overline{\U}\setminus\U$ of the open subset
$\U\subset\Stab(X)$ introduced in the last section.
The results of Section \ref{wall} show that $\partial\U$ is contained in
a locally-finite union of codimension one, real submanifolds of
$\Stab(X)$. A point
of the boundary will be called general if it lies on only one
of these submanifolds.

\begin{thm}
\label{hi}
Suppose $\sigma=(Z,\P)\in \partial\U$
is a general point of the boundary of $\U$. Then exactly one of the
following possibilities holds.
\begin{itemize}
\item[$(A^+)$] There is a rank $r$ spherical vector bundle
$A$ such that the only stable factors of the objects $\{\OO_x:x\in X\}$
in the stability condition $\sigma$ are $A$ and $T_A(\OO_x)$.
Thus the Jordan-H{\"o}lder filtration of each
$\OO_x$ is given by
\[ 0\lra A^{\oplus r}\lra \OO_x\lra T_A(\OO_x)\lra 0.\]

\item[$(A^-)$] There is a rank $r$ spherical vector bundle
$A$ such that the only stable factors of the objects $\{\OO_x:x\in X\}$
in the stability condition $\sigma$ are $A[2]$ and $T_A^{-1}(\OO_x)$.
Thus the Jordan-H{\"o}lder filtration of each
$\OO_x$ is given by
\[ 0\lra T_A^{-1}(\OO_x) \lra \OO_x\lra A^{\oplus r}[2]\lra 0.\]

\item[$(C_k)$] There is a non-singular rational curve $C\subset X$
and an integer $k$ such that $\OO_x$ is stable in
the stability condition $\sigma$ for $x\notin C$ and
such that the Jordan-H{\"o}lder filtration of $\OO_x$ for $x\in C$ is
\[ 0\lra \OO_C(k+1)\lra \OO_x\lra \OO_C(k)[1]\lra 0.\]
\end{itemize}
Moreover, a stability condition $\sigma$ satisfies $A^-$ precisely if the stability condition $T_A^{\,2}(\sigma)$ satisfies $A^+$,
 and similarly a stability condition $\sigma$ satisfies $C_k$ precisely if the stability condition $T_{\OO_C(k)}(\sigma)$ satisfies $C_{k-1}$. 
\end{thm}

\begin{pf}
If $\sigma\in\partial\U$ then each skyscraper sheaf $\OO_x$ is semistable of the same phase. Applying an element of $\grp$ it is enough to consider the case when
this phase is 1. Write $W(X)$ for the subset of $U(X)$ consisting of stability conditions for which each object $\OO_x$ is stable of phase 1.  Since $U(X)$ is invariant under the action of $\grp$ it follows that $\sigma$ lies in the closure of $W(X)$. 

\Step1
Since $\sigma\notin U(X)$ there must be some $x\in X$ such that
$\OO_x$ is semistable but not stable. Applying Lemma \ref{burb} below shows that there is a stable, spherical object $A\in\P(1)$ with a nonzero map
either $A\to\OO_x$ or $\OO_x\to A$.

Since $\sigma$ is general, the boundary of $\U$ at $\sigma$
is a codimension one
submanifold of $\Stab(X)$ at $\sigma$, which is given by $Z(A)/Z(\OO_x)\in\R_{>0}$.
Thus if $E$ is any object of $\D(X)$ one can
move a little bit along the boundary of $\U$
and assume that $Z(E)\not\in\R$, unless of course
the Mukai vector
$v(E)$ is a linear combination (over $\R$) of the vectors $v(A)$ and $v(\OO_x)$.

Let $L\subset\N(X)$ be the rank 2 sublattice spanned by $v(A)$ and $v(\OO_x)$. The fact that $A$ is spherical implies
that this lattice is primitive: if $v\in\N(X)$ lies in the linear span of the vectors $v(A)$ and $v(\OO_x)$ over $\R$ then
it also lies in the integral span, and hence is an element of $L$. Indeed one can write $v(A)=(r,\Delta,s)$ and $v(\OO_x)=(0,0,1)$.
If $v=\lambda v(A) + \mu v(\OO_x)$ lies in $\N(X)$ then $\lambda r\in\Z$, $\lambda\Delta\in\NS(X)$  and $\lambda s-\mu\in\Z$.
Clearly $\lambda$ and $\mu$ are rational. Since $\Delta^2-2rs=-2$ the elements $\Delta\in\NS(X)$ and $r\in\Z$ have no common divisor,
 so it follows that $\lambda$ and hence $\mu$ are integral. 

\Step2
Assume first that there is a nonzero map $A\to\OO_x$. 
Since $\sigma$ lies in the closure of $W(X)$, and $v(A)$ is primitive, by Proposition \ref{open} there are points $\tau=(W,\Q)\in W(X)$
arbitrarily close to $\sigma$ such that the object $A$ is also stable in
$\tau$. Thus one can take $f(\sigma,\tau)<\epsilon$ for some small $\epsilon>0$. The
nonzero map $A\to\OO_x$ then implies that $A\in\Q(\phi)$ with $1-\epsilon <\phi<1$ and it
follows from Lemma \ref{lemma}(a)
that $A$ fits into a triangle
\begin{equation*}
\xymatrix@C=.5em{ C \ar[rrrr] &&&& A  \ar[dll] \\
 && D \ar@{-->}[ull]} \end{equation*}
 with
$C=H^{-1}(A)[1]$ and $D=H^0(A)$.

Suppose for a contradiction that $C$ and $D$ are both nonzero.
Lemma \ref{lemma}(c) implies that $D\in\Q([-1,1])$ for all stability conditions $(W,\Q)\in W(X)$, and hence
$D\in\P([-1,-1])$ also. Similarly $C\in\P([0,2])$.
Lemma \ref{shipley} below shows that $C\in\P(0)$ and $D\in\P(1)$.
But now one could replace $\sigma$ by another nearby point of the
boundary of $\U$ and repeat the argument. Since $Z(C)$ and $Z(D)$ always
have to lie on
the real axis, and $\sigma$ was assumed to be general,
it follows that the Mukai vectors of $C$ and $D$ must
lie in the rank two sublattice $L$, so we can write
\[v(C)=\lambda_C v(A) + \mu_C v(\OO_x),\quad v(D)=\lambda_D v(A) + \mu_D v(\OO_x)\]
for integers $\lambda_C,\lambda_D,\mu_C,\mu_D\in\Z$,
with $\lambda_C+\lambda_D=1$ and $\mu_C+\mu_D=0$.

Moving again to the nearby stability condition $\tau=(W,\Q)\in W(X)$ considered above,
Lemma \ref{lemma}(c) shows that $C\in\Q((0,\epsilon))$ and
$D\in\Q((1-\epsilon,1])$. Also $A\in\Q(\phi)$ with  $1-\epsilon<\phi<1$  and of course $\OO_x\in\Q(1)$.
Thus the imaginary part of $W$ is positive on $C$ and $A$, non-negative on $D$ and zero on $\OO_x$. 
It follows that $\lambda_C>0$ and $\lambda_D\geq 0$. By Lemma \ref{rigid} the sheaf $D$ is rigid and so $v(D)^2<0$.
Since $v(\OO_x)^2=0$ it follows that $\lambda_D>0$ also which gives a contradiction.
Thus one of $C$ or $D$ is zero.

\Step3
If $D=0$ then $A=H^{-1}(A)[1]$ is concentrated in degree $-1$ and so
there could not be a nonzero map $A\to\OO_x$. Thus $C=0$ and hence $A$ is a sheaf.
If $A$ has positive rank then there is a short exact sequence of sheaves
\[0\lra T\lra A\lra Q\lra 0.\]
with $T$ torsion and $Q$ torsion-free.
Since $\sigma$ lies in the closure of $W(X)$, Lemma \ref{lemma}(c) shows that $T\in\P(\geq 0)$ and $Q\in\P([-1,1])$.
Lemma \ref{shipley} then shows that $T\in\P(0)$.
Again, since $\sigma$ is general it follows that $v(T)$ lies in the sublattice $L$. But $r(A)>0$ so this means that $v(T)$ is a multiple of $v(\OO_x)$
and so $T$ is supported in dimension 0. This too is impossible because all zero-dimensional sheaves lie in $\P(1)$.

Thus $A$ is either torsion or torsion-free.

\Step4
Suppose that $A$ is torsion-free with a map $A\to\OO_x$.
A result of Mukai \cite[Proposition 2.14]{Mu2} shows that $A$ must be
locally-free. Thus $\Hom_X(A,\OO_x)=\C^r$ for all $x\in X$ where $r$ is the rank of $A$.
Given $x\in X$ there is a short exact sequence in $\P(1)$ of the form
\[0\lra A'\lra \OO_x\lra B\lra 0\]
where $A'$ has all Jordan-H\"older factors isomorphic to $A$ and $\Hom_X(A,B)=0$.
Since $A$ is spherical (and hence has no self-extensions) it follows that $A'=A^{\oplus p}$ for some $p$,
and considering the long exact sequence obtained by applying the functor $\Hom_X(A,-)$ shows that $p=r$ and that the map
$A^{\oplus r}\to\OO_x$ is the canonical evaluation map.
It follows that $B=T_A(\OO_x)$.

Suppose for a contradiction that $B$ is not stable in $\sigma$.
Lemma \ref{burb} shows that
there is a stable spherical object $C\in\P(1)$ with either a nonzero
map $C\to B$ or a nonzero map $B\to C$.
The fact that $\sigma$ is general implies that $v(C)$
lies in the lattice $L$. Since $v(B)^2=-2=v(A)^2$ and $v(\OO_x)^2=0$, a quick calculation shows
that the only possibilities are $v(C)=\pm v(A)$, and since $A,C\in\P(1)$ we must take the positive sign.
The Riemann-Roch
theorem then
gives
$\eu(A,C)=2$, so that by Serre duality
there is either a map $C\to A$ or a map $A\to C$. Since $A$ and
$C$ are supposed to be stable of the same phase one concludes that
$C=A$. But
\[\Hom_X(A,B)=\Hom_X(T_A^{-1}(A),\OO_x)=\Hom_X(A[1],\OO_x)=0\]
and similarly $\Hom_X(B,A)=0$ which gives a contradiction. So $B$ is stable, and
the alternative $(A^+)$ applies.

\Step5
Suppose instead that $A$ is a torsion sheaf.
The Riemann-Roch theorem shows that $A$ is supported on a $(-2)-$curve
$C\subset X$. If
\[0\lra C\lra A\lra D\lra 0\]
is a short exact sequence of sheaves then Lemma \ref{lemma}(c) and Lemma \ref{shipley}
show that $C\in\P(0)$ and $D\in\P(1)$. Since $\sigma$ is general $v(C)$ and $v(D)$ are in the lattice $L$.
In particular it follows from this that $C$ must be non-singular, since if $C'\subset C$ is an irreducible component there is a surjection
of sheaves $A\to A|_{C'}$ which would lead to a non-trivial exact sequence as above. Riemann-Roch then shows that since $A$ is spherical one has
$A=\OO_C(k+1)$ for some integer $k$.
For any point $x\in C$
there is a unique map $A\to \OO_x$ with cone $B=\OO_C(k)[1]$. Since $A$ is
stable, this cone must also lie in $\P(1)$.

Suppose for a contradiction that $B$ is not stable. If $D\in\P(1)$  is a stable factor then by Lemma \ref{rigid} $D$ is spherical and since $\sigma$ is general
one concludes as in Step 4 that $v(D)$ lies in the lattice $L$ and hence $v(D)=\pm(0,C,l)$ for some $l$. Note that $Z(A)$ and $Z(B)$ lie on the negative real axis with
$Z(A)+Z(B)=-1$. Also $Z(D)$ must lie on the negative real axis with $|Z(D)|<|Z(B)|$. The only possibility is that $v(D)=v(A)$ which implies $D=A$. But if all stable factors of $D$ are isomorphic to $A$ one has $v(D)=nv(A)$ which is clearly false. Thus $\OO_C(k)[1]$ must be stable in $\sigma$ and
alternative $(C_k)$ applies.

\Step6
The cases when $A\in\P(1)$ with a map $\OO_x\to A$ are
dealt with by a similar analysis.
For points $\tau=(W,\Q)\in W(X)$
sufficiently close to $\sigma$ the object $A$ is also stable in
$\tau$, and the
nonzero map $\OO_x\to A$ then implies that $A\in\Q(\phi)$ with $1<\phi<1+\epsilon$.
Lemma \ref{lemma} then shows that there is a triangle
\begin{equation*}
\xymatrix@C=.5em{ C \ar[rrrr] &&&& A  \ar[dll] \\
 && D \ar@{-->}[ull]} \end{equation*}
 with
$C=H^{-2}(A)[2]$ and $D=H^{-1}(A)[1]$, and that moreover the sheaf $H^{-2}(A)$ is torsion-free.
Applying the argument of Step 2 one concludes that one of $C$ or $D$ is zero.

If $D=0$ then $A[-2]$ is a torsion-free sheaf and hence by Mukai's result \cite[Proposition 2.14]{Mu2} locally-free. The argument of Step 4 then shows that one is
in the situation $(A^-)$. If $C=0$ then the argument of Step 4 shows that $A[-1]$ is either torsion or torsion-free,
and the second possibility cannot hold because then $A[-1]$ would be locally-free and there could not be a nonzero map $\OO_x\to A$.
Thus $A[-1]$ is a torsion sheaf and the argument of Step 5 then shows that situation $(C_k)$ holds for some $k$.

\Step7
For the last statement, note that the triangles defining $T_A(\OO_x)$ and $T_A^{-1}(\OO_x)$ exist abstractly, so that a stability condition $\sigma$ satisfies $(A^-)$ precisely if the objects $A[2]$ and $T_A^{-1}(\OO_x)$ are stable of the same phase for all $x\in X$. Similarly $\sigma$ satisfies $(A^+)$ precisely if $A$ and $T_A(\OO_x)$ are stable of the same phase for all $x\in  X$. These two possibilities are clearly related by the equivalence $T_A^{2}$.
Similarly for the cases $(C_k)$ where one notes that the short exact sequence
of sheaves on $\PP^1$
\[0\lra\OO_{\PP^1}(k-1)\lra
\Hom_{\PP^1}(\OO_{\PP^1}(k),\OO_{\PP^1}(k+1))
\tensor\OO_{\PP^1}(k)\lra \OO_{\PP^1}(k+1)\lra 0\]
implies that $T_{\OO_C(k)}(\OO_C(k+1))=\OO_C(k-1)[1]$.
\end{pf}

The proof used the following three simple results.

\begin{lemma}
\label{burb}
Suppose
$\sigma=(Z,\P)\in\Stab(X)$ is a stability condition on $X$ and $0\neq E\in P(1)$
is a semistable object of phase 1.
If $\Hom^1_X(E,E)=0$ then all stable factors of $E$ are spherical objects. If $\Hom^1_X(E,E)= \C^2$ and $E$ is not stable
then there is a stable, spherical object $A\in\P(1)$ such that either $\Hom_X(A,E)$ or $\Hom_X(E,A)$ is nonzero.
\end{lemma}

\begin{pf}
Consider Jordan-H{\"o}lder filtrations
of $E$ in the finite length category $\P(1)$. Take a stable object $F\in\P(1)$
with a nonzero map $F\to E$. Grouping factors together, there is a short exact sequence
\[0\lra F'\to E\to G\lra 0\]
in $\P(1)$ such that all the Jordan-H{\"o}lder factors of $F'$ are
isomorphic to $F$, and $\Hom_X(F',G)=0$. Applying Lemma \ref{mukai}
shows that
\[\dim_{\C} \Hom^1_{X}(F',F')+\dim_{\C} \Hom^1_{X}(G,G)\leq \dim_{\C}\Hom^1(E,E)\leq 2,\]
and since both spaces have even dimension at least one of them is zero.
In the case $\Hom^1_X(E,E)=0$ both spaces are zero, so $v(F)^2<0$ and hence $F$ is spherical. Replacing $E$ by $G$ and repeating the argument gives the result.

Suppose instead that  $\dim_{\C}\Hom^1_X(E,E)= 2$. If $\Hom^1_X(F',F')=0$ then again $F$ is spherical and setting $A=F$ gives the result.
If $\Hom^1_X(G,G)=0$ then all the stable factors of $G$ are spherical, and so there is a stable, spherical object $A\in\P(1)$ with a
nonzero map $G\to A$, and hence a nonzero map $E\to A$.
\end{pf}

\begin{lemma}
\label{shipley}
Suppose $\sigma=(Z,\P)$ is a stability condition and $A\in\P(1)$ is stable. Given a triangle
\begin{equation*}
\xymatrix@C=.5em{ C \ar[rrrr] &&&& A  \ar[dll]^f \\
 && D \ar@{-->}[ull]} \end{equation*} 
with $C\in\P(\geq 0)$ and $D\in\P(\leq 1)$ and $f\neq 0$ then $C\in\P(0)$ and $D\in\P(1)$.
\end{lemma}

\begin{pf}
Considering the triangle $A\to D\to C[1]$ shows that $D\in\P(\geq 1)$ and hence $D\in\P(1)$. Since $A$ is simple in the abelian category
$\P(1)$ the map $f$ is a monomorphism in $\P(1)$ so the quotient $C[1]$ also lies in $\P(1)$.
\end{pf}

\begin{lemma}
\label{rigid}
If $E\in\D(X)$ satisfies $\Hom^1_X(E,E)=0$ then the same is true of each of its cohomology sheaves $H^i(E)$.
\end{lemma}

\begin{pf}
There is a
spectral sequence
\[E^{p,q}_2=\bigoplus_i\Ext^p_X(H^i(E),H^{i+q}(E))\implies
\Hom_X^{p+q}(E,E)\]
whose $E^{1,0}_2$ term survives to infinity.
\end{pf}

% ***************************************************************************
% ***************************************************************************
% ***************************************************************************
% ***************************************************************************

\section{Proof of the main theorem}
\label{proof}

It is now possible to use the results of the last two sections to prove
Theorem \ref{first}. The crucial point is that the images of
the closure of the set $\U$ under the elements of the group $\Aut\D(X)$ cover
the entire connected component $\Stab^\dagger(X)$.

\begin{lemma}
\label{nikita}
If $\sigma\in\Stab^\dagger(X)$ is a stability condition satisfying one of the assumptions $A^+$, $A^-$ or $C_k$ of Theorem \ref{hi}
then $\sigma$ lies in the boundary of $U(X)$.
\end{lemma}

\begin{pf}
Consider the $A^+$ case, the rest being similar. Set $B=T_A(\OO_x)$. Then  $A$ and $B$ are stable in $\sigma$ of the same phase $\phi$ and
there are short exact sequences
\[0\lra A^{\oplus r}\lra \OO_x\lra B\lra 0\]
in $\P(\phi)$. Applying an element of $\grp$ one can assume that $\phi=1$.

Since $\Stab^\dagger(X)$ is a good component and the objects $\{\OO_x:x\in X\}$ have bounded mass, one can take a wall and chamber decomposition as in Proposition \ref{lo}.
Clearly $\sigma$ must lie on a wall since all the objects $\OO_x$ are semistable but not stable.
There exists at least one chamber $\CC$ such that $\sigma$ lies in the closure of $\CC$ and such that $\Im W(A)/W(B)<0$ for stability conditions $(W,\Q)\in\CC$ close to $\sigma$. It will be enough to show  that all $\OO_x$ are stable in $\CC$ so that $\CC\subset\U$ and hence $\sigma$ lies in the closure of $\U$.

Consider a stability condition $\tau=(W,\Q)\in\CC$ such that $f(\sigma,\tau)<\ei$.  Let $\A$ be the abelian subcategory $\P((\tfrac{1}{2},\tfrac{3}{2}])\subset\D(X)$. If $\OO_x$ is  not stable there is a short exact sequence in $\A$
\[0\lra C\lra \OO_x\lra D\lra 0\]
with $\Im W(C)/W(D)>0$. One cannot have $\Im Z(C)/Z(D)>0$ because $\OO_x$ is semistable in $\sigma$. By the chamber structure
one must have $Z(C)/Z(D)\in\R_{>0}$ so that $C,D\in\P(1)$ and the above sequence is a short exact sequence in $\P(1)$. Considering stable factors of $\OO_x$ and noting that $\Hom_X(\OO_x,A)=0$
shows that all stable factors of $C$ are equal to $A$ so that $v(A)=nv(C)$. But this contradicts $\Im W(A)/W(B)<0$.
\end{pf}

\begin{prop}
\label{mendem}
The connected component $\Stab^\dagger(X)\subset\Stab(X)$ is mapped by $\Pi$
onto the open subset $\Pol^+_0(X)\subset\N(X)\tensor\C$.
\end{prop}

\begin{pf}
The fact that $\Pi(\Stab^\dagger(X))$ contains $\Pol^+_0(X)$ follows from the
fact that $\Pi$ is a covering map over $\Pol^+_0(X)$. The hard part is
to prove the reverse inclusion.

Take a stability condition $\sigma\in\Stab^\dagger(X)$. There is a
continuous path $\gamma:[0,1]\to\Stab^\dagger(X)$ such that $\gamma(0)\in\U$
and $\sigma=\gamma(1)$. One can find a compact
subset $\co\subset\Stab^\dagger(X)$ such that $\gamma([0,1])$ lies in its
interior. Let $S$ be the set of objects $E$ of $\D(X)$ such that $v(E)=v(\OO_x)$ and $E$ is semistable for
some stability condition $\sigma\in \co$. Then $m_{\sigma}(E)=|Z(E)|\leq m_{\sigma}(\OO_x)$ and since $\co$ is compact it follows that $S$ has bounded mass in the component $\Stab^\dagger(X)$. Consider the corresponding wall and chamber structure as in Section \ref{wall}.

Since $\pi$ is a local homeomorphism, and the walls are locally-finite, one can deform $\gamma$ a little so that
there are real numbers
$0=t_0<t_1<\cdots<t_n=1$ such that each interval
$I_i=(t_i,t_{i+1})$ is mapped by $\gamma$ into one of the chambers,
and such that each point $\gamma(t_i)$ for $0\leq i<n$ lies on only one wall.
Thus $\gamma(I_0)\subset U(X)$ and $\gamma(t_1)$ is a general point of the boundary of $U(X)$.

It follows from Theorem \ref{hi} and Lemma \ref{nikita} that at each general point $\sigma$ of the
boundary of $\U$ there is an autoequivalence $\Phi$ such that
$\Phi(\sigma)$ also lies in the boundary of $\U$. Moreover these
autoequivalences all preserve the class of $\OO_x$ in $\N(X)$.
The important point
 is that these autoequivalences reverse the orientation of
the boundary of $\U$.

Thus if $\sigma\in\partial \U$ is general of type $(A^-)$, then
locally near $\sigma$ the boundary of $\U$ is given by
the real quadric $\Im Z(A)/Z(\OO_x)=0$, and $\U$ is the
side where $\Im Z(A)/Z(\OO_x)>0$. Applying the equivalence $T_A^{\,2}$
gives a new stability condition on the boundary of $\U$, this time of
type $(A^+)$. Locally at
$T_A^{\,2}(\sigma)$ the boundary is still given by the equation
$\Im Z(A)/Z(\OO_x)=0$,
but now $\U$ is on the
side where $\Im Z(A)/Z(\OO_x)<0$.

On the other hand, if $\sigma\in\partial \U$
is general of type $(C_k)$, then the
boundary of $\U$ is given locally by $\Im Z(\OO_C)/Z(\OO_x)=0$,
with $\U$ being the side where $\Im Z(\OO_C)/Z(\OO_x)>0$. Applying
$T_{\OO_C(k)}$ identifies the $(C_k)$ part of the boundary of $\U$
with the $(C_{k-1})$ part, and the
fact that $T_{\OO_C(k)}$ acts on $\N(X)$ via a reflection shows
that again this identification is orientation-reversing.

Thus
there is an autoequivalence
$\Phi_1\in\Aut\D(X)$, preserving the class of $\OO_x$ in $\N(X)$, such
that $\Phi_1(\gamma(t_1))$ lies in  the boundary of $\U$, and for points $t>t_1$ close to $t$ one has
$\Phi_1(\gamma(t))\in \U$, which is to say $\Phi_1^{-1}(\OO_x)$ is stable in $\gamma(t)$.  By the chamber structure it follows that $\Phi_1(I_1)\subset \U$. Repeating the argument shows
that there is an autoequivalence $\Phi\in\Aut\D(X)$ preserving the
class of $\OO_x$ such that $\Phi(\sigma)$ lies in the closure of $\U$.

By Corollary \ref{liz}, the fact that $\Phi(\sigma)$ lies in the closure of $\U$ implies that
the real and imaginary parts of $\pi(\sigma)\in\N(X)\tensor\C$ spans a non-negative two-plane in
$\N(X)\tensor\R$.
Since this holds for any stability condition
in $\Stab^\dagger(X)$,
and the map $\Pi$ is open on a full component of $\Stab(X)$, it follows that
$\Pi(\Stab^\dagger(X))\subset\Pol(X)$. Since $\Stab^\dagger(X)$ is connected the image must in fact lie in $\Pol^+(X)$.

Finally, suppose that $\sigma=(Z,\P)\in\Stab^\dagger(X)$
satisfies $Z(\delta)=0$ for some
class $\delta\in\Delta(X)$. By the above one can assume that $\sigma$ lies in the boundary of $\U$.
According to Proposition \ref{lo} the boundary of $\U$ in a neighbourhood of $\sigma$ is made up of a finite union of codimension one submanifolds of $\Stab(X)$ each passing through $\sigma$. By Theorem \ref{hi}, each of these components is of the form $Z(A)/Z(\OO_x)\in \R_{>0}$ for some spherical object $A$. Moreover each $A$ is stable of the same phase as $\OO_x$ at a general point of the corresponding boundary component, and hence semistable at $\sigma$.

Fix one component $\Gamma$ of the boundary of $\U$ near $\sigma$ and suppose $\delta$ lies in the corresponding  rank two sublattice $L\subset \N(X)$ spanned by $v(A)$ and $v(\OO_x)$. If $r(A)>0$ then since $v(A)$ and $\delta$ both lie in $\Delta(X)$, one has $v(A)=\pm \delta$, which is impossible because $A$ is semistable at $\sigma$ whereas $Z(\delta)=0$. If $r(A)=0$ then the component $\Gamma$ is of the form $(C_k)$ for some $(-2)-$curve $C\subset X$ and some integer $k$, so that at a general point the Jordan-H{\"o}lder filtration of $\OO_x$ for points $x\in C$ are
of the form 
\[ 0\lra \OO_C(k+1)\lra \OO_x\lra \OO_C(k)[1]\lra 0.\]
Then $\OO_C(k)$ and $\OO_C(k+1)$ are at least semistable at $\sigma$ so have non-vanishing central charge, and hence
$0<|Z(\OO_C(k)|<|Z(\OO_x)|$. 
But if $\delta\in L$ then $\delta=\pm v(\OO_C(l))$ for some integer $l$ and then  
$Z(\delta)=0$ implies that $Z(\OO_C(k))=(k-l)Z(\OO_x)$ which gives a contradiction.

Since the image of $\Stab^\dagger(X)$ is open in $\N(X)\tensor \C$ it now follows that there exist stability conditions in $\Stab^\dagger(X)$ arbitrarily close to $\sigma$ which are contained in $\U$ and which still satisfy $Z(\delta)=0$. But this is impossible by Corollary \ref{liz}. 
\end{pf}

The following is a
restatement of Theorem \ref{first}.

\begin{thm}
The map $\Pi\colon\Stab^\dagger(X)\to\Pol^+_0(X)$ is a covering map.
The subgroup of $\Aut^0\D(X)$ fixing the connected component $\Stab^\dagger(X)\subset\Stab(X)$
acts freely on $\Stab^\dagger(X)$ and is the group of deck transformations.
\end{thm}

\begin{pf}
The fact that $\Pi$ is a covering map is Proposition \ref{cover}. It
is clear that the
given subgroup of $\Aut^0\D(X)$ acts on $\Stab^\dagger(X)$ preserving the map
$\Pi$. What we are required to show is that given stability conditioins $\sigma$ and $\tau$
in $\Stab^\dagger(X)$ with the same central charge there is a unique $\Phi\in\Aut^0\D(X)$ with $\Phi(\tau)=\sigma$.
Since $\pi$ is a covering map, it is enough to check this for a fixed $\sigma\in\Stab^\dagger(X)$ which we may as well assume lies in $\U$.
 
To prove uniqueness suppose $\Phi\in\Aut^0\D(X)$ is such that $\Phi(\sigma)=\sigma$.
Thus if $E$ is stable in  $\sigma$ of a given phase then so is $\Phi(E)$.
By Lemma  \ref{lemma}, the only objects
$E\in\D(X)$ which are stable in $\sigma$ with Mukai vector $v(E)=v(\OO_x)$ are the skyscraper sheaves
$\OO_x$ themselves together with their even shifts.
It follows that $\Phi$ takes skyscrapers to
skyscrapers which implies that $\Phi(E)=f^*(E\tensor L)$ for some
$L\in\Pic(X)$ and some $f\in\Aut(X)$. Since $\Phi\in\Aut^0\D(X)$ it
follows (using the Torelli theorem) that $\Phi$ is the identity.

Now assume that $\sigma$ and $\tau$ have the same central charge and $\sigma\in \U$.
By the argument of Proposition \ref{mendem}, there is an
autoequivalence $\Phi\in\Aut\D(X)$ such that
$\Phi(\tau)$ lies in the closure of $\U$, so moving $\sigma$
a bit one can assume that $\sigma$ and $\tau$ both lie in $\U$.
 Moreover one can assume that $\Phi$ preserves the
class of $\OO_x$ in $\N(X)$. Composing $\Phi$ with a twist by a line
bundle, which preserves the set $\U$, one can assume that $\Phi$ also
preserves the class of $\OO_X$. Thus the action of $\Phi$ on $H^*(X,\Z)$
preserves the decomposition
\[H^*(X,\Z)=H^0(X,\Z)\oplus H^2(X,\Z)\oplus H^4(X,\Z).\]
Since $\sigma$ and $\Phi(\tau)$ both
lie in $\U$, the induced Hodge isometry of $H^2(X,\Z)$ is effective
\cite[Proposition VIII.3.10]{BPV}.
It follows from the Torelli theorem \cite[Theorem VIII.11.1]{BPV} that
composing $\Phi$ with an automorphism of $X$, which again
preserves $\U$, one can assume that $\Phi$ acts trivially on
$H^*(X,\Z)$. But now $\sigma$ and $\Phi(\tau)$ have the same central
charge and both lie in $\U$. It follows from Corollary \ref{liz} that composing $\Phi$ with an even shift one has $\Phi(\tau)=\sigma$
which completes the proof.
\end{pf}

% ***************************************************************************
% ***************************************************************************
% ***************************************************************************
% ***************************************************************************

\section{The large volume limit}

The aim of this section is to study the class of stable objects
in the stability condition $\sigma\in V(X)$
corresponding to a point $\exp(\beta+i\omega)\in\Pol^+_0(X)$
in the limit as $\omega\to\infty$.
This is what physicists would refer to as a large
volume limit, and string theory predicts that the BPS branes in this
limit are just Gieseker stable sheaves. In fact, in the prescence of a
nonzero B-field, Gieseker stability gets twisted. This leads to a
notion of stability first introduced by Matsuki and Wentworth
\cite{MW}.

Throughout this section $\beta,\omega\in\NS(X)\tensor\R$ will be a
fixed pair of $\R$-divisor classes with $\omega\in\Amp(X)$ ample.
Given a torsion-free
sheaf $E$ on $X$ with Mukai vector $v(E)=(r(E),\cl_1(E),s(E))$
define
\[\mu_{\beta,\omega}(E)
=\frac{\big(c_1(E)-r(E)\beta\big)\cdot \omega}{r(E)}\quad
\text{ and }\quad\nu_{\beta,\omega}(E)=
\frac{s(E)-\cl_1(E)\cdot \beta}{r(E)}.\]
Note that $\mu_{\beta,\omega}(E)=\mu_\omega(E)-\beta\cdot\omega$.
The following definition reduces to Gieseker
stability in the case $\beta=0$.

\begin{defn}
A torsion-free sheaf $E$ on $X$ is said to be
twisted semistable with respect to the pair $(\beta,\omega)$
if
\[\mu_{\beta,\omega}(A)<\mu_{\beta,\omega}(E)\text{ or }\big(\mu_{\beta,\omega}(A)=\mu_{\beta,\omega}(E)
\text{ and }\nu_{\beta,\omega}(A)\leq\nu_{\beta,\omega}(E)\big)\]
for all subsheaves $0\neq A\subset E$.
\end{defn}

Note that a twisted semistable sheaf is, in particular, slope
semistable.
The mathematical reason for introducing twisted stability
is that unlike slope stability, Gieseker stability is not preserved by
twisting by line bundles. Thus if $\beta\in\NS(X)$ is the first
Chern class of a line bundle $L$, then a torsion-free sheaf $E$ is
Gieseker semistable with respect to $\omega$ if and only if $E\tensor
L$ is
twisted semistable with respect to the pair $(\beta,\omega)$. The above
definition just generalises this idea to arbitrary elements $\beta\in\NS(X)\tensor\R$.

\begin{prop}
Fix a pair $\beta,\omega\in\NS(X)\tensor\QQ$ with $\omega\in\Amp(X)$
ample.
For
integers $n\gg 0$ there is a unique stability condition
$\sigma_n\in\U$ satisfying  $\Pi(\sigma_n)=\exp(\beta+in\omega)$. 
Suppose $E\in\D(X)$ satisfies
\[r(E)> 0\ \text{ and }\ 
(\cl_1(E)-r(E)\beta)\cdot \omega>0.\]
Then $E$ is semistable
in $\sigma_n$ for all
$n\gg 0$ precisely if
$E$ is a shift of a $(\beta,\omega)$-twisted semistable sheaf on $X$.
\end{prop}

\begin{pf}
Providing $(n\omega)^2>2$ one has
$\exp(\beta+in\omega)\in\KKo(X)$ so it follows from
Proposition \ref{below} that
there is a unique stability condition
$\sigma_n\in\U$ satisfying  $\Pi(\sigma_n)=\exp(\beta+in\omega)$.
Note that each of the stability conditions $\sigma_n=(Z_n,\P_n)$ has the same
heart $\A(\beta,\omega)=\P_n((0,1])$.
Note also that if $E$ is a nonzero sheaf on $X$ then equation
$(\star)$ of Section \ref{con} shows that
\begin{equation*}
\lim_{n\to\infty}\frac{1}{\pi}\arg Z_n(E)=\begin{cases}
0 &\text{ if }\supp(E)=X, \\
\ha &\text{ if }\dim\supp(E)=1, \\
1 &\text{ if }\supp(E)=0.
\end{cases}
\end{equation*}

Take an object $E\in\D(X)$ with $r(E)>0$ and
$\mu_{\beta,\omega}(E)>0$.
First suppose that $E$ is semistable in
$\sigma_n$ for all $n\gg 0$. Applying a shift one may as well assume
that $E$ lies in $\A(\beta,\omega)$.
According to Lemma \ref{lemma}, $E$ has
non-vanishing cohomology sheaves in just two degrees, and there is a
short exact sequence in $\A(\beta,\omega)$
\[ 0\lra H^{-1}(E)[1]\lra E\lra H^0(E)\lra 0.\]
Now $H^{-1}(E)$ is a torsion-free sheaf, so according to the asymptotic
formula
above, the phase of $E$ tends to 0 in the limit $n\to\infty$,
whereas, if the object $H^{-1}(E)[1]$ is nonzero,
its phase must tend
to 1. Since $E$ is semistable in
$\sigma_n$ for $n\gg 0$ it follows
that $H^{-1}(E)=0$ and hence $E$ is a sheaf.
A similar argument with the asymptotic formula shows that $E$ must be
torsion-free.
Note that the $\mu_{\beta,\omega}$-semistable factors of $E$ all have positive slope.

Suppose $E$ is not $(\beta,\omega)$-twisted semistable. Then
there is a destabilising sequence
\[0\lra A\lra E\lra B\lra 0\]
of sheaves on $X$ such that $A$ and $B$ lie in $\A(\beta,\omega)$ and $A$ is a $\mu_{\omega}$-semistable sheaf with
$\mu_{\beta,\omega}(A)>\mu_{\beta,\omega}(E)$.
Rewriting equation $(\star)$ of Section \ref{con} gives
\[\frac{Z_n(E)}{r(E)}-\frac{Z_n(A)}{r(A)}=
-\big(\nu_{\beta,\omega}(E)-\nu_{\beta,\omega}(A)\big)
+in\big(\mu_{\beta,\omega}(E)-\mu_{\beta,\omega}(A)\big).\]
Since the phases of $A$ and $E$ tend to zero, this implies
$\arg Z_n(A)>\arg Z_n(E)$ for all $n\gg 0$ which contradicts semistability of $E$.
Thus $E$ is $(\beta,\omega)$-twisted semistable.

For the converse, suppose
that $E$ is a $(\beta,\omega)$-twisted
semistable torsion-free sheaf
with $\mu_{\beta,\omega}(E)>0$.
In particular, $E$ is $\mu_{\omega}$-semistable so that
$E\in\A(\beta,\omega)$.
Suppose
\[0\lra A\lra E\lra B\lra 0\]
is a short exact sequence in $\A(\beta,\omega)$.
Taking cohomology gives a
long exact sequence of sheaves
\begin{equation*}
\tag{$\dagger$}
0\lra H^{-1}(B)\lra A\lra E\lra H^0(B)\lra 0.
\end{equation*}

By definition of the category $\A(\beta,\omega)$ the sheaf $H^{-1}(B)$
has slope $\mu_{\beta,\omega}\leq 0$ whereas $A$ has slope
$\mu_{\beta,\omega}> 0$.
Since $E$ is semistable it follows that
$\mu_{\beta,\omega}(A)\leq\mu_{\beta,\omega}(E)$ and if
equality holds then $\nu_{\beta,\omega}(A)\leq\nu_{\beta,\omega}(E)$.
Lemma \ref{bel} shows that the set of possible values of
$\nu_{\beta,\omega}(A)$ is
bounded above. The phase of $E$ tends to zero as $n\to\infty$ and
moreover the real part of $Z_n(E)$ tends to $+\infty$.
It follows that $\arg Z_n(A)\leq
\arg Z_n(E)$
for all subobjects $A\subset E$ in $\A(\beta,\omega)$ and all $n\gg 0$.
Thus $E$ is semistable in
$\sigma_n$ for $n\gg 0$.
\end{pf}

\begin{lemma}
\label{bel}
Let $E$ be a $(\beta,\omega)$-twisted semistable
torsion-free
sheaf with
$\mu_{\beta,\omega}(E)>0$. Then the set of values of $\nu_{\beta,\omega}(A)$
as $A$ ranges through all
nonzero subobjects of $E$ in $\A(\beta,\omega)$ is bounded above.
\end{lemma}

\begin{pf}
Since $\beta$ is rational, Theorem \ref{yosh} shows
that there exist torsion-free $\mu_{\omega}$-semistable sheaves $P$ with
Mukai vector $(r,r\beta,s)$ for some $s\in\Z$.
Decomposing into
stable factors and taking a double dual one may assume that $P$ is in fact
$\mu_{\omega}$-stable and locally-free.

The Riemann-Roch formula gives
\[\eu(P,A)=s(P)r(A)+r(P)\big( s(A)- c_1(A)\cdot\beta\big).\]
Comparing with the formula for $\nu_{\beta,\omega}(A)$ shows that  
it is enough to give an upper bound for $\eu(P,A)/r(A)$.
Consider the exact sequence $(\dagger)$ and put
$D=H^{-1}(B)$.
Since $A$
has Harder-Narasimhan factors of positive slope $\mu_{\beta,\omega}$
there can be
no maps $A\to P$. Thus it suffices to bound the quotient $\dim_{\C}
\Hom_X(P,D)/r(D)$.

Since the Harder-Narasimhan factors of $D$ have non-positive slope
$\mu_{\beta,\omega}$, 
any nonzero map $f\colon P\to D$ is an injection with torsion-free quotient.
Indeed, if $f$ is not injective, then since $P$ is $\mu_{\beta,\omega}$-stable, the image of $f$
would have
strictly positive slope $\mu_{\beta,\omega}$ which is impossible. Thus there is a short exact sequence
\[0\lra P\lra D\lra Q\lra 0.\]
If $Q$ has torsion subsheaf $T$, then the induced map $D\to Q/T$ has kernel $K$ fitting into an exact sequence
\[0\lra P\lra K\lra T\lra 0.\]
If $T$ is supported in dimension one then $K$ has strictly positive slope $\mu_{\beta,\omega}$ which again is impossible. But if $T$ is supported in dimension zero then the above sequence splits which is impossible since $D$ is torsion-free. Hence $T=0$.

Replacing $D$ by $Q$ and proceeding by induction on $r(D)$ it is easy to see that
\[\dim_\C\Hom_X(P,D)\leq r(D)/r(P). \]
which completes the proof.
\end{pf}

% ***************************************************************************
% ***************************************************************************
% ***************************************************************************
% ***************************************************************************

\section{Stability conditions on abelian surfaces}

Suppose that instead of a K3 surface one considers an abelian surface
$X$. The theory developed in this paper carries over virtually
unchanged to this case, because the only features of a K3 surface
which were used were homological properties such as the Riemann-Roch
theorem or Serre duality. The most important difference between
abelian surfaces and K3 surfaces from the point of view of this paper
is the following

\begin{lemma}
\label{tat}
If $X$ is an abelian surface then $\D(X)$ has no spherical objects.
\end{lemma}

\begin{pf}
In fact there are no nonzero objects $E\in\D(X)$ with $\Hom^1_X(E,E)=0$.
The basic reason is that there can be no rigid objects on a
torus because of the continuous automorphism group.

By Lemma \ref{rigid} it is enough to consider the case when $E$ is a sheaf.
Note that if $E$ is a vector bundle, then $\OO_X$ is a
direct summand of $E^{\vee}\tensor E$, so that 
\[\C^2=H^1(X,\OO_X)\subset H^1(X,E^{\vee}\tensor E)=\Ext^1_X(E,E).\]
For a general sheaf $E$ one can twist by a sufficiently ample line
bundle and apply the Fourier-Mukai transform \cite{Mu1} to obtain a
vector bundle with the same $\Ext$-algebra.
\end{pf}

Suppose then that $X$ is an abelian surface over $\C$.
Consider the
even cohomology lattice
\[H^{2*}(X,\Z)=H^0(X,\Z)\oplus H^2(X,\Z)\oplus H^4(X,\Z),\]
equipped with the Mukai bilinear form.
This is an even and non-degenerate
lattice of signature $(4,4)$.
The Todd class of $X$ is trivial,
so the Mukai vector of an object $E\in\D(X)$ is the triple
\[v(E)=\ch(E)=(r(E),c_1(E), \ch_2(E))\in \N(X)=\Z\oplus\NS(X)\oplus\Z
.\]
As before, $\N(X)$ has signature $(2,\rho)$ and
$\Pol^+(X)\subset\N(X)\tensor\C$ is defined to be one component of the
set of vectors $\mho\in\N(X)\tensor\C$ which span positive definite
two-planes.

It is a simple consequence of Lemma \ref{tat} that, up to a shift, every
autoequivalence $\Phi\in\Aut\D(X)$ takes skyscraper sheaves to sheaves
\cite[Corollary 2.10]{BM}. It follows
immediately
that $\Aut^0\D(X)$ is generated by the
double shift $[2]$, together with twists by elements
of $\Pic^0(X)$, and pull-backs by automorphisms of $X$ acting trivially on $H^*(X,\Z)$.
It also follows that no element of $\Aut\D(X)$ exchanges the
two components of $\Pol(X)$.

As before one can define an open subset $\U\subset\Stab(X)$.
This time however, because of the lack of spherical objects,
$\U=\Stab^\dagger(X)$ is actually a connected component of
$\Stab(X)$. Thus up to the action of $\grp$, every stability condition in
$\Stab^\dagger(X)$ is obtained from the construction of Section \ref{con}.

Note that twists by elements of $\Pic^0(X)$ and pull-backs by automorphisms of $X$
acting trivially on $H^*(X,\Z)$ take
skyscrapers to skyscrapers, 
so these elements of $\Aut^0\D(X)$ act trivially on
$\Stab^\dagger(X)$.
Note also that $\Pol^+(X)$ is a $\GL(2,\R)$-bundle over the contractible space
\[\{\beta+i\omega\in \NS(X)\tensor \C: \omega^2>0\}\]
so that $\pi_1 \Pol^+(X)=\Z$. Putting these observations together gives

\begin{thm}
Let $X$ be an abelian surface over $\C$.
There is a connected component $\Stab^\dagger(X)\subset\Stab(X)$
which is mapped by $\Pi$ onto the open subset
$\Pol^+(X)\subset\N(X)\tensor\C$.
Moreover, the induced map
\[\Pi\colon \Stab^\dagger(X)\lra\Pol^+(X)\]
is the universal cover, and the group of deck transformations
is generated by the double shift functor $[2]$.
\end{thm}

One can also show directly that the action of the group $\Aut \D(X)$ on
$\Stab(X)$ preserves
the connected component $\Stab^\dagger(X)$. Thus the analogue of Conjecture
\ref{conj} holds in the abelian surface case.

% ***************************************************************************
% ***************************************************************************
% ***************************************************************************
% ***************************************************************************

\section{Final remarks}

So far we have not mentioned moduli spaces of stable objects. One
might hope that something of the following sort is true.

\begin{conj}
Given a nonsingular projective variety $X$,
a stability condition $\sigma\in\Stab(X)$,  and a
numerical equivalence class $\alpha\in\N(X)$,
there exists a coarse moduli
space $\M_{\sigma}(\alpha)$
of objects in $\D(X)$ of type $\alpha$ which are semistable in
$\sigma$.
\end{conj}

Note that as it stands this conjecture is rather vague
since
it does not
specify what category the coarse moduli space should live in.

Let us now suppose for definiteness that $X$ is a K3 surface.
The group
$\Aut \D(X)$ of exact autoequivalences
$\Phi\colon \D(X)\to\D(X)$
acts on the space $\Stab(X)$ in such a way that an object
$E\in\D(X)$ is stable in a given stability condition $\sigma$ if and
only if $\Phi(E)$ is stable in the stability condition $\Phi(\sigma)$.
The point I wish to make here is that
this action suggests the idea that a 
Fourier-Mukai transform should be thought of
as being analagous to a change in polarisation.

Recall the standard technique for using Fourier-Mukai transforms to
compute moduli spaces of stable bundles, going back to Mukai's original
work on the subject \cite{Mu3}.
One considers a universal family of Gieseker stable
bundles $\{E_s:s\in S\}$ on $X$
of some fixed numerical type $\alpha$ and applies a Fourier-Mukai
transform $\Phi$ to obtain a new family of objects $\Phi(E_s)$ of $\D(X)$ 
of some different numerical type $\beta$. In certain special examples
the objects $\Phi(E_s)$ are also Gieseker stable bundles, thus
giving an isomorphism of moduli spaces $\M(\alpha)\to\M(\beta)$.
In general of course, the objects $\Phi(E_s)$ are just complexes and
this approach then fails. But from the perspective of this paper, one
could say that the
objects $\Phi(\E_s)$ are in fact always stable, but only
with respect to some transformed
stability condition on $\D(X)$.

In fact, as we saw in Proposition \ref{lo},
for a given numerical class
$\alpha$, the space $\Stab(X)$ splits into a collection of
walls and chambers, and an object $E\in\D(X)$ of type $\alpha$
can only become stable or unstable by crossing from one
chamber to another. Thus for general $\sigma\in\Stab(X)$
the space $\M_{\sigma}(\alpha)=\M_{\CC}(\alpha)$
only depends on the
chamber $\CC$ in which $\sigma$ lies.
Furthermore, if $\sigma$ lies sufficiently close to the large
volume limit point of $\Stab(X)$, an object of $\D(X)$ of type
$\alpha$ or $\beta$ is stable in $\sigma$
precisely if it is a shift of a Gieseker stable
sheaf. Thus the functor
$\Phi$ gives an identification between the moduli space of
Gieseker stable sheaves $\M(\alpha)$ and the moduli space
$\M_{\Phi(\sigma)}(\beta)$, and
 the question of whether $\Phi$ preserves Gieseker
stability can be addressed by considering whether $\sigma$ and $\Phi(\sigma)$
lie in the same chamber.

Understanding the relationship between
$\M_{\Phi(\sigma)}(\beta)$ and $\M(\beta)$ directly seems to be rather
difficult. From our point of view this
is because of the
complicated geometry of the wall and chamber structure on $\Stab(X)$.
But assuming $\Stab(X)$ is connected, one can choose
a sequence
of adjacent chambers $\CC_1,\cdots, \CC_N$ for $\beta$
with $\sigma\in \CC_1$ and $\Phi(\sigma)\in \CC_N$. It is then
tempting
to speculate that by analysing
wall-crossing phenomena one might hope to prove that
in good cases, in analogy to \cite{MW}, there is a
birational equivalence 
\[\M(\alpha)=\M_{\CC_N}(\beta)\dashrightarrow\cdots\dashrightarrow
\M_{\CC_2}(\beta)\dashrightarrow
\M_{\CC_1}(\beta)=\M(\beta)\]
arising as a sequence of  flops.

% ***************************************************************************
% ***************************************************************************
% ***************************************************************************
% ***************************************************************************

\bigskip

\noindent Department of Mathematics,
University of Sheffield,
Hicks Building, Hounsfield Road, Sheffield, S3 7RH, UK.

%\smallskip

\noindent email: {\tt t.bridgeland@sheffield.ac.uk}


\begin{thebibliography}{101}

%\bibitem{AP} D. Abramovich and A. Polishchuk, Sheaves of t-structures
%and valuative criteria for stable complexes, preprint (2003).

%\bibitem{AM} P. Aspinwall and D. Morrison, String theory on K3
%surfaces, Mirror symmetry, II, 703--716, 
%AMS/IP Stud. Adv. Math., 1, 
%Amer. Math. Soc., Providence, RI, 1997, also hep-th/9404151.

\bibitem{BPV} W. Barth, C. Peters and A. Van de Ven,
Compact complex surfaces. Ergebnisse der Mathematik
(3), 4. Springer-Verlag, Berlin, 1984.

\bibitem{BBD} A. Beilinson, J. Bernstein and P. Deligne, Faisceaux
Pervers, Ast{\'e}risque {100}, Soc. Math de France (1983).

\bibitem{BM} T. Bridgeland and A. Macioica, Complex surfaces with equivalent derived categories,
Math. Z. 236 (2001).

\bibitem{B2} {T. Bridgeland and A. Maciocia}, Fourier-Mukai transforms
for $K3$ and elliptic fibrations. J. Algebraic Geom. 11
 (2002), no. 4.


\bibitem{B} {T. Bridgeland,} Stability conditions on triangulated categories,
to appear in Ann. Math. (2007) Preprint math.AG/0212237.

\bibitem{B3} T. Bridgeland, Spaces of stability conditions, ********* to appear.



\bibitem{Do4} M.R. Douglas,  Dirichlet branes, homological mirror symmetry, and stability. Proceedings of the International Congress of Mathematicians, Vol. III (Beijing, 2002), 395--408, Higher Ed. Press, Beijing, 2002, also  math.AG/0207021.

\bibitem{GM} S.I. Gelfand and Yu. I. Manin, Methods of Homological
Algebra, Springer-Verlag (1996).

\bibitem{HRS} D. Happel, I. Reiten, S.O. Smal{\o}, Tilting in abelian
categories and quasitilted algebras. Mem. Amer. Math. Soc. 120 (1996),
no. 575.

%\bibitem{In} M. Inaba, Towards a definition of moduli of complexes of
%coherent sheaves on a projective scheme, J. Math.
%Kyoto Univ. 42 (2002), no. 2, 317--329.

\bibitem{MW} K. Matsuki, R. Wentworth, Mumford-Thaddeus principle on
the moduli space of vector bundles on an algebraic surface.
Internat. J. Math. 8 (1997), no. 1, 97--148.

\bibitem{Mu1} S. Mukai, Duality between $\D(X)$ and $\D(\hat{X})$ with its
application to Picard sheaves, Nagoya Math. J. {\bf 81} (1981) 153--175.

\bibitem{Mu2} { S. Mukai,} On the moduli space of bundles on K3 surfaces I,
in: Vector Bundles on Algebraic Varieties, M.F. Atiyah et al., Oxford
University Press (1987), 341-413.

\bibitem{Mu3} S. Mukai, Fourier functor and its application to the
moduli of bundles on an abelian variety. Algebraic geometry, Sendai,
1985, 515--550, Adv. Stud.
Pure Math., 10 (1987).

\bibitem{Or} { D.O. Orlov,} Equivalences of derived categories and K3
surfaces, J. Math. Sci. (New York) { 84} 5 (1997) 1361-1381,
also alg-geom 9606006.

\bibitem{ST} P. Seidel, R.P. Thomas, Braid group actions on derived
categories of coherent sheaves, Duke Math. J. 108 (2001), no. 1,
37--108.

%\bibitem{Si} C.T. Simpson, Moduli of representations of the fundamental
%group of a smooth projective variety I, Publ. Math. I.H.E.S. {79}
%(1994) 47-129.

\bibitem{Sz} B. Szendr\H oi, Diffeomorphisms and families of
Fourier--Mukai transforms in mirror symmetry,
Applications of algebraic geometry to coding theory, physics and
computation, 317--337, NATO Sci. Ser., Kluwer 2001.


\bibitem{Yo2} K. Yoshioka, Chamber structure of polarizations and
the moduli of stable sheaves on a ruled surface, Internat.
J. Math. 7 (1996), no. 3, 411--431, also math.AG/9409008.

\bibitem{Yo1} K. Yoshioka, Moduli spaces of stable sheaves on abelian
surfaces, Math. Ann. 321 (2001), no. 4, 817--884.

\end{thebibliography}
\end{document}